\documentclass[10pt]{article}

\usepackage{amsmath}
\usepackage{amsfonts}
\usepackage{amssymb}
\usepackage[english]{babel}
\usepackage{amsthm}
\usepackage{extarrows}
\usepackage{amscd}
\usepackage{tikz}
\usepackage{mathdots}
\usetikzlibrary{arrows,shapes}

\def\pf{\par\noindent {\bf Proof}~\par\noindent}
\def\qed{~\hfill{$\square$}\pagebreak[1]\par\medskip\par}

\newcommand{\mR}{\mathbb{R}}
\newcommand{\mC}{\mathbb{C}}

\newcommand{\mZ}{\mathbb{Z}}
\newcommand{\mE}{\mathbb{E}}
\newcommand{\mS}{\mathbb{S}}
\newcommand{\mH}{\mathbb{H}}
\newcommand{\mI}{\mathbb{I}}
\newcommand{\mJ}{\mathbb{J}}
\newcommand{\mK}{\mathbb{K}}
\newcommand{\mQ}{\mathbb{Q}}

\newcommand{\mV}{\mathbb{V}}
\newcommand{\mM}{\mathbb{M}}
\newcommand{\mU}{\mathbb{U}}
\newcommand{\mW}{\mathbb{W}}

\newcommand{\mcG}{\mathcal{G}}
\newcommand{\mcP}{\mathcal{P}}
\newcommand{\mcQ}{\mathcal{Q}}
\newcommand{\mcM}{\mathcal{M}}
\newcommand{\mcH}{\mathcal{H}}
\newcommand{\mcI}{\mathcal{I}}

\newcommand{\mcS}{\mathcal{S}}
\newcommand{\mcT}{\mathcal{T}}
\newcommand{\mcE}{\mathcal{E}}

\newcommand{\gf}{\mathfrak{f}}
\newcommand{\gfd}{\mathfrak{f}^{\dagger}}
\newcommand{\gsl}{\mathfrak{sl}}
\newcommand{\gsp}{\mathfrak{sp}}
\newcommand{\gso}{\mathfrak{so}}
\newcommand{\gosp}{\mathfrak{osp}}

\newcommand{\gH}{\mathfrak{H}}
\newcommand{\geg}{\mathfrak{g}}
\newcommand{\gl}{\mathfrak{gl}}

\newcommand{\ol}{\overline}
\newcommand{\olz}{\ol{z}}

\newcommand{\uX}{\underline{X}}

\newcommand{\uz}{\underline{z}}
\newcommand{\uzJ}{\underline{z}^J}
\newcommand{\uzJd}{\underline{z}^{\dagger J}}
\newcommand{\uzd}{\underline{z}^{\dagger}}

\newcommand{\upz}{\partial_{\uz}}
\newcommand{\upzJ}{\partial_{\uz}^J}

\newcommand{\upzd}{\partial_{\uz}^{\dagger}}
\newcommand{\upzJd}{\partial_{\uz}^{\dagger J}}

\newcommand{\nz}{\vert\underline{z}\vert}

\newcommand{\p}{\partial}
\newcommand{\dirac}{\underline{\p}}

\newcommand{\SH}{\mathcal{H}^S}
\newcommand{\SHD}{\mathcal{H}^{\dagger S}}
\newcommand{\SE}{\mcE}
\newcommand{\SED}{\mcE^{\dagger}}

\newcommand{\onehalf}{\frac{1}{2}}
\newcommand{\onequarter}{\frac{1}{4}}

\newtheorem{theorem}{Theorem}
\newtheorem{lemma}{Lemma}
\newtheorem{proposition}{Proposition}
\newtheorem{definition}{Definition}
\newtheorem{remark}{Remark}
\newtheorem{corollary}{Corollary}

\begin{document}

\title{Fischer Decomposition for $\gosp(4|2)$--monogenics in Quaternionic Clifford Analysis}
\author{F.\ Brackx$^\ast$, H.\ De Schepper$^\ast$, D.\ Eelbode$^{\ast\ast}$, R.\ L\'{a}vi\v{c}ka$^\ddagger$ \& V.\ Sou\v{c}ek$^\ddagger$}

\date{\small{$^\ast$ Clifford Research Group, Dept. of Math. Analysis, 
Faculty of Engineering and Architecture,\\ Ghent University,
Building S22, Galglaan 2, B--9000 Gent, Belgium\\
$\ast\ast$ University of Antwerp, Middelheimlaan 2, Antwerpen, Belgium\\
$^\ddagger$ Charles University in Prague, Faculty of Mathematics and Physics, Mathematical Institute\\
Sokolovsk\'a 83, 186 75 Praha, Czech Republic}}

\maketitle


\begin{abstract}
Spaces of spinor--valued homogeneous polynomials, and in particular spaces of spinor--valued spherical harmonics, are decomposed in terms of irreducible representations of the symplectic group Sp$(p)$. These Fischer decompositions involve spaces of homogeneous, so--called $\gosp(4|2)$--monogenic polynomials, the Lie super algebra $\gosp(4|2)$ being the Howe dual partner to the symplectic group Sp$(p)$. In order to obtain Sp$(p)$--irreducibility this new concept of $\gosp(4|2)$--monogenicity has to be introduced as a refinement of quaternionic monogenicity; it is defined by means of the four quaternionic Dirac operators, a scalar Euler operator $\SE$ underlying the notion of symplectic harmonicity and a  multiplicative Clifford algebra operator $P$ underlying the decomposition of spinor space into symplectic cells. These operators $\SE$ and $P$, and their hermitian conjugates, arise naturally when constructing the Howe dual pair $\gosp(4|2) \times {\rm Sp}(p)$, the action of which will make the Fischer decomposition multiplicity free.

\end{abstract}


\section{Introduction}
\label{introduction}


In a similar way as hermitian Clifford analysis in Euclidean space $\mR^{2n}$ of even dimension arises as a refinement of Euclidean Clifford analysis by introducing a complex structure on  $\mR^{2n}$, quaternionic Clifford analysis originates as a further refinement by the introduction of a so--called hypercomplex structure on  Euclidean space $\mR^{4p}$, the dimension now being a fourfold. A hypercomplex structure $\mQ$ (see Section \ref{QCA}) consists of three, freely chosen but fixed, mutually anti--commuting complex structures $\mI_{4p}$, $\mJ_{4p}$ and $\mK_{4p}$, i.e.\ elements of the SO$(4p)$ group, the composition of which mimicks the multiplication rules for the quaternionic units $i, j, k$ in the quaternion algebra $\mH$:
$$
\mI_{4p}^2 = - {\bf 1}, \ \mJ_{4p}^2 = - {\bf 1}, \ \mK_{4p}^2 = - {\bf 1}, \ \mI_{4p} \mJ_{4p} = - \mJ_{4p} \mI_{4p} = \mK_{4p}
$$
The subgroup SO$_{\mQ}(4p)$ of SO$(4p)$, consisting of those matrices which commute with the complex structures in $\mQ$, is isomorphic with the symplectic group Sp$(p)$. So--called Stein--Weiss projections of the gradient operator (see e.g.\ \cite{stein}) yield four differential operators $\upz$, $\upzd$, $\upzJ$, $\upzJd$, which are invariant under the Sp$(p)$-action. Their simultaneous null solutions are called {\em quaternionic monogenic} functions; they have recently been studied in e.g.\ \cite{PSS, eel2, DES}. For the fundaments of this function theory and a precise description of the involved systems  of differential equations we refer to \cite{paper1,paper2}.\\[-2mm]

The aim of the present paper is to establish a Fischer decomposition of spaces of spinor--valued homogeneous polynomials, and in particular spaces of spinor--valued spherical harmonics, in terms of Sp$( p)$--irreducibles (Section \ref{fischer}). To that end the value space, i.e.\ spinor space $\mS$, has to be decomposed into so--called symplectic cells which are fundamental representations of Sp$(p)$. For a detailed description of this splitting we refer to \cite{paper1}; in Section \ref{spinordecomposition}  the construction is briefly recalled.  However, spaces of quaternionic monogenic functions with values in a symplectic cell are not Sp$( p)$--irreducible, whence the need for (two) additional operators, a scalar differential operator $\mcE$ (or its hermitian conjugate $\SED$) and a Clifford multiplication operator $P$ (or its hermitian conjugate $Q$), to refine the notion of quaternionic monogenicity to so--called $\gosp(4|2)$--monogenicity (Section \ref{osp42monogenic}), which reflects by its name the fact that the Howe dual partner to the symplectic group Sp$(p)$ precisely is the Lie superalgebra $\gosp(4|2)$. In fact, as opposed to the traditional development where function theory comes first and the Fischer decomposition and the related Howe dual pair later legitimize the choices made in defining and studying the systems of equations, the definition of $\gosp(4|2)$--monogenicity by means of the six operators $\upz$, $\upzd$, $\upzJ$, $\upzJd$, $\mcE$ (or $\SED$) and $P$ (or $Q$),  is dictated by this Lie superalgebra. To make the paper self--contained the basics of Clifford algebra and of Euclidean and hermitian Clifford analysis are recalled in Sections \ref{cliffordalgebra} and \ref{EHCA} respectively.


\section{Preliminaries on Clifford algebra}
\label{cliffordalgebra}


For a detailed description of the structure of Clifford algebras we refer to e.g.\ \cite{port}. Here we only recall the necessary basic notions. The real Clifford algebra $\mathbb{R}_{0,m}$ is constructed over the vector space $\mathbb{R}^{0,m}$ endowed with a non--degenerate quadratic form of signature $(0,m)$, and generated by the orthonormal basis $(e_1,\ldots,e_m)$. The non--commutative multiplication in $\mathbb{R}_{0,m}$ is governed by the rules 
\begin{equation}\label{multirules}
e_{\alpha} e_{\beta} + e_{\beta} e_{\alpha} = -2 \delta_{\alpha \beta} \ \ , \ \ \alpha,\beta = 1,\ldots ,m
\end{equation}
As a basis for $\mathbb{R}_{0,m}$ one takes for any set $A=\lbrace j_1,\ldots,j_h \rbrace \subset \lbrace 1,\ldots,m \rbrace$ the element $e_A = e_{j_1} \ldots e_{j_h}$, with $1\leq j_1<j_2<\cdots < j_h \leq m$, together with $e_{\emptyset}=1$, the identity element; the dimension of $\mR_{0,m}$ is $2^m$. Any Clifford number $a$ in $\mathbb{R}_{0,m}$ may thus be written as $a = \sum_{A} e_A a_A$, $a_A \in \mathbb{R}$, or still as $a = \sum_{k=0}^m \lbrack a \rbrack_k$, where $\lbrack a \rbrack_k = \sum_{|A|=k} e_A a_A$ is the so--called $k$--vector part of $a$. Real numbers correspond to the zero--vector part of the Clifford numbers, while Euclidean space $\mathbb{R}^{0,m}$ is embedded in $\mathbb{R}_{0,m}$ by identifying $(X_1,\ldots,X_m)$ with the Clifford $1$--vector $\uX = \sum_{\alpha=1}^m e_{\alpha}\, X_{\alpha}$, for which it holds that $\uX^2 = - |\uX|^2 = - r^2$.\\[-2mm]

When allowing for complex constants, the generators $(e_1,\ldots, e_{m})$, still satisfying (\ref{multirules}), produce the complex Clifford algebra $\mathbb{C}_{m} = \mathbb{R}_{0,m} \oplus i\, \mathbb{R}_{0,m}$. Any complex Clifford number $\lambda \in \mathbb{C}_{m}$ may thus be written as $\lambda = a + i b$, $a,b \in \mathbb{R}_{0,m}$, this form leading to the definition of the hermitian conjugation $\lambda^{\dagger} = (a +i b)^{\dagger} = \overline{a} - i \overline{b}$, where the bar notation stands for the Clifford conjugation in $\mathbb{R}_{0,m}$, i.e.\ the main anti--involution for which $\overline{e}_{\alpha} = -e_{\alpha}$, $\alpha=1, \ldots,m$. The hermitian conjugation leads to a hermitian inner product on $\mathbb{C}_{m}$ given by $(\lambda,\mu) = \lbrack \lambda^{\dagger} \mu \rbrack_0$, and its associated norm $|\lambda| = \sqrt{ \lbrack \lambda^{\dagger} \lambda \rbrack_0} = ( \sum_A |\lambda_A|^2 )^{1/2}$.\\[-2mm]


\section{Euclidean and hermitian Clifford analysis: the basics}
\label{EHCA}


The central notion in standard Clifford analysis  is that of a {\em monogenic function}. This is a continuously differentiable function defined in an open region of Euclidean space $\mR^m$, taking values in the Clifford algebra $\mR_{0,m}$, or subspaces thereof, and vanishing under the action of the Dirac operator $\dirac = \sum_{\alpha=1}^m e_{\alpha}\, \p_{X_{\alpha}}$, a vector--valued first order differential operator, which can be seen as the Fourier or Fischer dual of the Clifford variable $\uX$. Monogenicity is the higher dimensional counterpart of holomorphy in the complex plane. As the Dirac operator factorizes the Laplacian: $\Delta_m = - \dirac^2$, Clifford analysis can also be regarded as a refinement of harmonic analysis. \\[-2mm]

It is important to note that the Dirac operator is invariant under the action of the $\mathrm{Spin}(m)$--group which doubly covers the $\mathrm{SO}(m)$--group, whence this framework is usually referred to as Euclidean (or orthogonal) Clifford analysis. A key result in this function theory is the following Fischer decomposition which fully reflects the importance of the underlying Spin$(m)$ symmetry. Note that the functions under consideration do not take their values in the whole Clifford algebra but in spinor space. Indeed, a Clifford algebra may be decomposed as a direct sum of isomorphic copies of a spinor space $\mS$, which, abstractly, may be defined as a minimal left ideal in the Clifford algebra and moreover is an irreducible Spin$(m)$ representation, see also Remark \ref{spinor-remark} below.

\begin{proposition}
\label{propmonogenicfischer}
The space $\mcP(\mR^m;\mS)$ of spinor--valued polynomials can be decomposed as
\begin{equation}
		\label{monogenicfischer}
\mcP(\mR^{m};\mS) = \bigoplus_{k=0}^\infty \bigoplus_{j=0}^{\infty} \uX^j \; \mcM_{k}(\mR^m;\mS)
\end{equation}
with $\mcM_k(\mR^m;\mS)$ the space of $k$--homogeneous monogenic  polynomials, where for each  $j=0,1, \ldots $ and each $k=0,1, \ldots $ the space $\uX^j \, \mcM_k(\mR^m;\mS)$ is an irreducible representation of $\mathrm{Spin}(m)$.
\end{proposition}

\noindent The Fischer decomposition (\ref{monogenicfischer}) however has a drawback: it is not multiplicity free since all spaces $\uX^j \, \mcM_k(\mR^m;\mS), j=0,1,\ldots$ are isomorphic 
Spin$(m)$-representations. Whence the need for the corresponding Howe dual pair, the action of which will lead to a multiplicity free Fischer decomposition. In this case (see \cite{howe}) the Howe dual partner is the Lie superalgebra
$\gosp(1|2) = \geg_0 \oplus \geg_1 = \gsl(2) \oplus \geg_1$
with
$$
\gsl(2) \cong \mathrm{Alg}_{\mC} \left( \mE + \frac{m}{2}, \onehalf r^2, - \onehalf \Delta_m \right)
$$
where $\mE = \sum_{\alpha=1}^{m} \, X_\alpha \p_{X_\alpha}$ is the Euler operator, and $\geg_1 \cong \mathrm{span}_{\mC} \left( \uX, \dirac  \right)$. The $\mZ$--grading of the orthosymplectic Lie superalgebra $\mathfrak{osp}(1|2)$ has the structure
$$
\mathfrak{osp}(1|2) = \mcG_{-2} \oplus  \mcG_{-1} \oplus \mcG_{0} \oplus \mcG_{1} \oplus \mcG_{2}
$$
where
$$
\mcG_{-2} = \mathrm{Alg}_{\mC} \left( r^2 \right)
$$
$$
\mcG_{-1} = \mathrm{span}_{\mC} \left( \uX \right)
$$
$$
\mcG_{0} = \mcH = \mathrm{Alg}_{\mC} \left( \mE  + \frac{m}{2} \right)
$$
$$
\mcG_{1} = \mathrm{span}_{\mC} \left( \dirac  \right)
$$
$$
\mcG_{2} = \mathrm{Alg}_{\mC} \left( \Delta_m \right)
$$
and it is easily verified that, as is well-known: $[g_{2j} , g_{k}] = g_{2j+k}$ and $\{g_{2j+1} , g_{2k+1}\} = g_{2(j+k+1)}$, for $g_i \in \mcG_i$, tacitly assuming the result to be zero whenever the sum of the indices is meaningless. By means of this Howe dual pair, Proposition \ref{propmonogenicfischer} may be reformulated as follows.
\begin{theorem}
\label{theomonogenicfischer}
Under the joint action of $\gosp(1|2) \times {\rm Spin}(m)$, the space $\mcP(\mR^{m};\mS)$ is isomorphic to the
multiplicity free irreducible direct sum
$$
\bigoplus_{k=0}^{\infty} \widetilde{\mI}_k \otimes \mM_k \ ,
$$
where $\widetilde{\mI}_k$ is the irreducible $\gosp(1|2)$--module realized by $\{\uX^j \, M_k : M_k \in \mcM_k(\mR^m; \mS), j=0,1,\ldots   \}$ and  $\mM_k$ denotes the irreducible {\rm Spin}$(m)$--module isomorphic to the space of $\mS$--valued $k$--homogeneous monogenic polynomials.
\end{theorem}
\noindent In view of Theorem \ref{theomonogenicfischer} we may now call a (standard) monogenic function also $\gosp(1|2)$--monogenic.\\[-2mm]

Taking the dimension of the underlying Euclidean vector space $\mR^m$ to be even: $m=2n$, renaming the variables: $(X_1,\ldots,X_{2n}) = (x_1,y_1,x_2,y_2,\ldots,x_n,y_n)$ and considering the standard complex structure $\mI_{2n}$, i.e. the complex linear real $\mbox{SO}(2n)$--matrix 
$$
\mI_{2n} = \mathrm{diag} \begin{pmatrix} \phantom{-} 0 & 1 \\ -1 & 0 \end{pmatrix}$$
for which $\mI_{2n}^2 = - E_{2n}$, $E_{2n}$ being the identity matrix, we define the rotated vector variable
$$
\uX_\mI = \mI_{2n}[\uX] = \mI_{2n} \left[ \sum_{k=1}^n (e_{2k-1} x_k + e_{2k} y_k ) \right] = \sum_{k=1}^n \mI_{2n} [e_{2k-1}] x_k + \mI_{2n} [e_{2k}] y_k = \sum_{k=1}^n (-y_k e_{2k-1} + x_k e_{2k})
$$
and, correspondingly, the rotated Dirac operator
$$
\dirac_\mI = \mI_{2n} [ \dirac] = \sum_{k=1}^n ( - \p_{y_k} e_{2k-1} + \p_{x_k} e_{2k})
$$

A differentiable function $F$ then is called {\em hermitian monogenic} in some region $\Omega$ of $\mR^{2n}$, if and only if in that region $F$ is a solution of the system
\begin{equation}
\dirac F = 0 = \dirac_\mI F
\label{hmon}
\end{equation}
Observe that this notion of hermitian monogenicity does not involve the use of complex numbers, but instead, could be developed as a real function theory. There however is an alternative approach making use of the projection operators $\frac{1}{2} ({\bf 1} \pm i \, \mI_{2n})$ and thus involving a complexification.
In this approach we consider in the complexification $\mC^{2n}$ of $\mR^{2n}$ the so--called Witt basis vectors, given by
$$
\gf_k = - \frac{1}{2} ({\bf 1} - i \, \mI_{2n}) [e_{2k-1}] \quad \mathrm{and} \quad \gfd_k =  \frac{1}{2} ({\bf 1} + i \, \mI_{2n}) [e_{2k-1}] \quad (k=1,\ldots,n)
$$
They satisfy the Grassmann identities
$$
\gf_j \gf_k + \gf_k \gf_j = 0 , \quad \gf_j^\dagger \gf_k^\dagger + \gf_k^\dagger \gf_j^\dagger = 0, \qquad j,k=1,\ldots,n
$$
including their isotropy $\gf_j^2 = (\gf_j^\dagger)^2 = 0$, $j=1,\ldots, n$, as well as the duality identities
$$
\gf_j \gf_k^\dagger + \gf_k^\dagger \gf_j = \delta_{jk}, \qquad j,k=1,\ldots,n
$$
The Witt basis vectors $(\gf_1,\ldots,\gf_{n})$ on the one hand, and $(\gf_1^\dagger,\ldots,\gf_{n}^\dagger)$ on the other, respectively span isotropic subspaces $W$ and $W^\dagger$ of $\mC^{2n}$, such that $\mC^{2n} = W \oplus W^\dagger$, those subspaces being eigenspaces of the complex structure $\mI_{2n}$ with respective eigenvalues $-i$ and $i$. They also generate the respective Grassmann algebras $\mC \Lambda_{n}$ and $\mC \Lambda_{n}^\dagger$. Note that the $\cdot^\dagger$--notation corresponds to the hermitian conjugation in the Clifford algebra $\mC_{2n}$.\\[-2mm]

We now consider the vector variables
\begin{eqnarray*}
\uz &=&  - \frac{1}{2} ({\bf 1} - i \, \mI_{2n} ) [\uX] \ = \ \sum_{k=1}^n ( x_k \gf_k + y_k (i \gf_k) ) = \sum_{k=1}^n (x_k + i y_k) \gf_k = \sum_{k=1}^n z_k \gf_k \\
\uzd &=& \phantom{-} \frac{1}{2} ({\bf 1} + i \, \mI_{2n} ) [\uX] \ = \ \sum_{k=1}^n ( x_k \gfd_k + y_k (-i \gfd_k) ) = \sum_{k=1}^n \overline{z}_k \gfd_k
\end{eqnarray*}
and, correspondingly, the hermitian Dirac operators
\begin{eqnarray*}
2 \, \upzd &=&  - \frac{1}{2} ( {\bf 1} - i \, \mI_{2n})[\dirac] = \sum_{k=1}^n ( \gf_k \p_{x_k} + i \gf_k \p_{y_k}) = \sum_{k=1}^n \gf_k ( \p_{x_k} + i \p_{y_k} ) = 2 \sum_{k=1}^n \p_{\overline{z}_k} \gf_k \\
2 \, \upz &=&  \phantom{-} \frac{1}{2} ( {\bf 1} + i \, \mI_{2n})[\dirac] = \sum_{k=1}^n ( \gfd_k \p_{x_k} - i \gfd_k \p_{y_k}) = \sum_{k=1}^n \gfd_k ( \p_{x_k} - i \p_{y_k} ) = 2 \sum_{k=1}^n \p_{z_k} \gfd_k
\end{eqnarray*}
As $2(\upz - \upzd) = \dirac$ and $2(\upz+\upzd) = i \, \dirac_\mI$, it follows that the system (\ref{hmon}) is equivalent with the system
\begin{equation}
\upz F = 0 = \upzd F
\label{hmoneq}
\end{equation}
\noindent Again it suffices to consider functions with values in a spinor space $\mS$, which now is an irreducible Spin$(2n)$ representation, and will be explicitly realized below. 

\begin{remark}\label{spinor-remark}
Here it is necessary to point out that there exist two kinds of (complex) spinors, called Dirac spinors and Weyl spinors. Their definition and behaviour as a module for the group Spin$(m)$, or its associated orthogonal Lie algebra $\gso(m)$, depends on the parity of the underlying dimension $m$. When the dimension $m$ is even, which is the case in the present hermitian framework since $m=2n$, and will also be the case in the quaternionic setting where $m=4p$, the Dirac spinor space $\mS$ is usually defined as $\mS := \mC_{m} \, I$, with $I \in \mC_{m}$ a suitable primitive idempotent, i.e.\ $I^2 = I$ (see below). Under the multiplicative action
$$ \textup{Spin}(m) \times \mS \rightarrow \mS : (s,\psi) \mapsto s\psi $$
this spinor space then decomposes into two irreducible representations $\mS^\pm$ for the spin group {\rm Spin}$(m)$, called the spaces of positive and negative Weyl spinors respectively, which can be defined as $\mS^+ := \mC_{m}^+ \, I$, with $\mC_{m}^+ \subset \mC_{m}$ the even subalgebra, and $\mS^- := \mC_{m}^- \, I$, with $\mC_{m}^- \subset \mC_{m}$ the odd subspace. It is indeed easily seen that the action of the spin group is well-defined on these Weyl spinor spaces, as {\rm Spin}$(m) \subset \mC_{m}^+$. It turns out that both spinor spaces define non-equivalent representations of the spin group with half-integer highest weights $(\frac{1}{2},\ldots,\frac{1}{2},\pm \frac{1}{2})$. These observations already apply to the even dimensional case of Euclidean Clifford analysis; in particular we should reconsider Proposition \ref{propmonogenicfischer} as well as Theorem \ref{theomonogenicfischer}, where, technically speaking, the space $\mcM_k(\mR^m,\mS)$ is not irreducible with respect to the spin action when $m$ is even. However, we will stick to the above formulation, where if necessary both sides of the decomposition  (\ref{monogenicfischer}) can be restricted to either positive or negative spinors. For the sake of convenience, also in what follows we will continue to work with Dirac spinors $\mS := \mS^+ \oplus \mS^-$, taking into account that the Dirac operator and its various `deformations' (which are also vector-valued differential operators) map positive spinors to negative ones and vice versa.
\end{remark}

With the self-adjoint idempotents 
$$
I_j = \gf_j \gf_j^\dagger = \frac{1}{2} ( 1- i e_{2j-1} e_{2j} ), \qquad j=1,\ldots,n
$$
we compose the primitive self--adjoint idempotent $I = I_1 I_2 \cdots I_{n}$, leading to the realization of the spinor space as $\mS = \mC_{2n} I $. Since $\gf_j I =0$, $j=1,\ldots,n$, we also have $\mS \simeq \mC \Lambda_{n}^\dagger I$. When decomposing the Grassmann algebra as 
$$
\mC \Lambda_{n}^\dagger = \bigoplus_{r=0}^{n} \left ( \mC \Lambda_{n}^\dagger \right )^{(r)}
$$ 
into its so--called homogeneous parts,
where $\left (\mC \Lambda_{n}^\dagger \right )^{(r)}$ is spanned by all products of $r$ Witt basis vectors out of $(\gf_1^\dagger,\ldots,\gf_{n}^\dagger)$, the spinor space $\mS$ accordingly decomposes into
$$
\mS = \bigoplus_{r=0}^{n} \mS^r, \qquad \mbox{with\ } \mS^r \simeq \left ( \mC \Lambda_{n}^\dagger \right)^{(r)} I
$$
These homogeneous parts $\mS^r$, $r=0,\ldots,n$, of spinor space provide models for fundamental $\mbox{U}(n)$--representations (see \cite{howe}) and for fundamental $\gsl_{n}(\mC)$-representations (see \cite{partI}, \cite{Ee}), whence functions with values in a homogeneous spinor subspace $\mS^r, r=0,\ldots,2n$ are to be considered. As to the group invariance of the operators $\upz$ and $\upzd$ defining hermitian monogenicity, they are invariant under the action of $\mbox{SO}_{\mI}(2n)$, i.e.\ the subgroup of $\mbox{SO}(2n)$ consisting of those matrices which commute with the complex structure $\mI_{2n}$. This subgroup $\mbox{SO}_\mI(2n)$ inherits a twofold covering by the subgroup $\mbox{Spin}_{\mI}(2n)$ of $\mbox{Spin}(2n)$, consisting of those elements of $\mbox{Spin}(2n)$ which commute with 
\begin{equation}
s_{\mI} = s_1 \ldots s_{n},  \quad \mbox{where \ } s_j = \frac{\sqrt{2}}{2} ( 1 + e_{2j-1} e_{2j} ),  \; j=1,\ldots,n
\label{sI}
\end{equation}
The element $s_{\mI}$ itself obviously belongs to $\mbox{Spin}_{\mI}(2n)$ and corresponds, under the double covering, to the complex structure $\mI_{2n} $. As $\mbox{SO}_{\mI}(2n)$ is isomorphic with the unitary group $\mbox{U}(n)$, whence up to this isomorphism $\mbox{Spin}_{\mI}(2n)$ provides a double cover of $\mbox{U}(n)$, we may say that $\mbox{U}(n)$ is the fundamental group underlying the function theory of hermitian monogenic functions. This U$(n)$--symmetry is apparent in the following result concerning the Fischer decomposition in terms of homogeneous hermitian monogenic polynomials (see \cite{howe}). To that end we introduce the space $\mcI$ of all Spin$_{\mI}(2n)$--invariant polynomials, which may be proved by invariance theory (see e.g. \cite{good,howe}) to be spanned by all words in the letters $\uz$ and $\uzd$: $ \mcI = \mbox{span}_{\mC} \left(1,\uz, \uzd, \uz\: \uzd, \uzd \uz, \uz\:\uzd\uz, \uzd\uz\:\uzd, \uz\:\uzd\uz\:\uzd, \uzd\uz\:\uzd\uz, \cdots \right)$, or alternatively
$$
\mcI = \mbox{span}_{\mC} \left( w_l^{(i)}(\uz,\uzd) : l = 0, 1, 2, \ldots , i = 1,2  \right)
$$
where $w_0^{(1)} = w_0^{(2)} = 1$ and
$$
\begin{array}{cccccccccl}
	w_{2j}^{(1)}(\uz,\uzd) & = & (\uz \uzd)^j & = & \nz^{2j-2}\uz\:\uzd & & & w_{2j+1}^{(1)}(\uz,\uzd) & = & \nz^{2j}\uz\\[2mm]
	w_{2j}^{(2)}(\uz,\uzd) & = & (\uzd \uz)^j & = & \nz^{2j-2}\uzd \uz & & & w_{2j+1}^{(2)}(\uz,\uzd) & = & \nz^{2j}\uzd\ .
\end{array}
$$

\begin{proposition}
\label{prophermitfischer}
The space $\mcP(\mR^{2n};\mS)$ of spinor--valued polynomials can be decomposed according to the $\mathrm{Spin}_\mI(2n)$ action as
\begin{equation}
\label{hermitfischer}
\mcP(\mR^{2n};\mS) =  \bigoplus_{a,b=0}^{\infty} \bigoplus_{r=0}^{n} \left( \mcM_{a,b}^{r}(\mR^{2n};\mS^r)  \oplus \bigoplus_{l=1}^{\infty} \ \bigoplus_{i=1,2} \ w_l^{(i)}(\uz,\uzd) \ \mcM_{a,b}^{r}(\mR^{2n};\mS^r) \right)
\end{equation}
with $\mcM_{a,b}^{r}(\mR^{2n};\mS^r)$ the space of $(a,b)$--homogeneous hermitian monogenic  polynomials  in the complex variables $(z_1, \cdots, z_n, \olz_1, \cdots, \olz_n)$ with values in the homogeneous spinor subspace $\mS^r$.
\end{proposition}

\noindent
Again this Fischer decomposition (\ref{hermitfischer}) in terms of irreducible $\mathrm{Spin}_\mI(2n)$--invariant subspaces is not multiplicity free, whence the need for the action of the Howe dual pair. It turns out that the dual partner to $\mathrm{U}(n)$ is the Lie super algebra
$
\mathfrak{sl}(1|2) = \geg_0 \oplus \geg_1 = \gl(2) \oplus \geg_1 = \gsl(2) \oplus \mC \oplus \geg_1
$
with
$$
\geg_1 \cong \mathrm{span}_{\mC} \left( \uz, \uzd, \upz, \upzd  \right)
$$
and
$$
\gsl(2) \cong \mathrm{Alg}_{\mC} \left( \mE_{\uz} +  \mE_{\uz}^\dagger  + n, \onehalf r^2, - \onehalf \Delta_{2n} \right),
\qquad
\mC \cong  \mathrm{span}_{\mC} \left( \mE_{\uz}^\dagger - \mE_{\uz} + n -2\beta \right)
$$
where $\mE_{\uz}$ and  $\mE_{\uz}^\dagger$ are the Euler operators corresponding to the hermitian variables $\uz$ and $\uzd$, respectively:
$$
\mE_{\uz} = \sum_{j=1}^n \, z_j \, \p_{z_j} \quad \mE_{\uz}^\dagger = \sum_{j=1}^n \, \olz_j \, \p_{\olz_j}
$$
which split the standard Euler operator as $\mE = \mE_{\uz} + \mE_{\uz}^\dagger$,  and $\beta$ is the so--called spin--Euler operator
$$
\beta = \sum_{j=1}^n \, \gfd_j \gf_j
$$
which measures the homogeneity degree of a homogeneous spinor subspace. The $\mZ$--gradation of this Lie superalgebra $\gsl(1|2)$ has the structure
$$
\gsl(1|2) = \mcG_{-2} \oplus  \mcG_{-1} \oplus \mcG_{0} \oplus \mcG_{1} \oplus \mcG_{2}
$$
where
$$
\mcG_{-2} = \mathrm{Alg}_{\mC} \left( \nz^2 \right)
$$
$$
\mcG_{-1} = \mathrm{span}_{\mC} \left( \uz, \uzd \right)
$$
$$
\mcG_{0} = \gH = \mathrm{Alg}_{\mC} \left(\mE_{\uz} +  \mE_{\uz}^\dagger  + n,  \mE_{\uz} - \mE_{\uz}^\dagger + n -2\beta \right)
$$
$$
\mcG_{1} = \mathrm{span}_{\mC} \left( \upz, \upzd  \right)
$$
$$
\mcG_{2} = \mathrm{Alg}_{\mC} \left( \Delta_{2n} \right)
$$
By means of this Howe dual pair, Proposition \ref{prophermitfischer} may be reformulated as follows.

\begin{theorem}
\label{theohermitfischer}
Under the joint action of $\gsl(1|2) \times {\rm Spin}_\mI(2n)$, the space $\mcP(\mR^{2n};\mS)$ of spinor--valued polynomials is isomorphic to the multiplicity free irreducible direct sum
$$
\bigoplus_{a,b=0}^{\infty}\bigoplus_{r=0}^{n} \widetilde{\mI}_{a,b,r} \otimes \mM_{a,b,r}
$$
where $\mM_{a,b,r}$ denotes the irreducible {\rm Spin}$_\mI(2n)$--module isomorphic to the space of $\mS^{r}$--valued $(a,b)$--homogeneous hermitian monogenic polynomials, and $\widetilde{\mI}_{a,b,r}$ denotes the $\gsl(1|2)$--irreducible module isomorphic to
$$
\widetilde{\mI}_{a,b,r} \cong \bigoplus_{w \in \mcI} \, {\rm span}_{\mC} \left( w M_{a,b}^r  \right)
$$
$M_{a,b}^r$ being an arbitrary but fixed polynomial in $\mcM_{a,b}^r(\mR^{2n};\mS^r)$.
\end{theorem}

\begin{remark}
There is an alternative expression, involving harmonic functions, for the the $\gsl(1|2)$--irreducible module $\widetilde{\mI}_{a,b,r}$, namely
$$
\widetilde{\mI}_{a,b,r} \cong \left( \mU^1 \oplus \mU^2 \right) \oplus \left( \mW^1 \oplus \mW^2  \right)
$$
with
\begin{eqnarray*}
\mU^1 &=& \mathrm{span}_\mC \{ |\uz|^{2l} \, M_{a,b}^r :  M_{a,b}^r \in \mcM_{a,b}^r(\mR^{2n};\mS^r), l = 0, 1,\ldots    \}\\
\mU^2 &=& \mathrm{span}_\mC \{ |\uz|^{2l} \, ((a+r)\uz\uzd - ((b+n-r)\uzd \uz) \, M_{a,b}^r :  M_{a,b}^r \in \mcM_{a,b}^r(\mR^{2n};\mS^r), l = 0, 1,\ldots    \}\\
\mW^1 &=& \mathrm{span}_\mC \{ w_{2j+1}^{(1)} \, M_{a,b}^r :  M_{a,b}^r \in \mcM_{a,b}^r(\mR^{2n};\mS^r), j = 0, 1,\ldots    \}\\
\mW^2 &=& \mathrm{span}_\mC \{ w_{2j+1}^{(2)} \, M_{a,b}^r :  M_{a,b}^r \in \mcM_{a,b}^r(\mR^{2n};\mS^r), j = 0, 1,\ldots    \}
\end{eqnarray*}
In fact it is this alternative expression which has been used in \cite{howe} to prove Theorem \ref{theohermitfischer}.
\end{remark}

In view of Theorem \ref{theohermitfischer} we may now call a hermitian monogenic function also $\gsl(1|2)$--monogenic.


\section{Quaternionic Clifford Analysis}
\label{QCA}


A refinement of hermitian Clifford analysis is obtained by considering the hypercomplex  structure $\mQ = (\mI_{4p},\mJ_{4p},\mK_{4p})$ on $\mR^{4p} \simeq \mC^{2p} \simeq \mH^p$, the dimension $m=2n=4p$ now assumed to be a 4-fold. To this end we introduce, next to the complex structure $\mI_{4p}$, a second one, $\mJ_{4p}$, given by
$$
\mJ_{4p} = \mbox{diag} \, \left ( \begin{array}{cccc} &  & 1 & \\ & & & -1 \\  -1 & & &  \\ & 1 & & \end{array} \right )
$$
Clearly $\mJ_{4p}$ belongs to $\mbox{SO}(4p)$, with $\mJ_{4p}^2 = -E_{4p}$, and it anti--commutes with $\mI_{4p}$. 
Then, quite naturally, there is a third $\mbox{SO}(4p)$--matrix $\mK_{4p} = \mI_{4p} \, \mJ_{4p} = - \mJ_{4p} \, \mI_{4p}$ for which $\mK_{4p}^2 = -E_{4p}$ and which anti--commutes with both $\mI_{4p}$ and $\mJ_{4p}$. It turns out that 
$$
\mK_{4p} = \mbox{diag} \, \left ( \begin{array}{cccc} &  &  & -1 \\ & & -1 &  \\   & 1 & &  \\  1 & & & \end{array} \right )
$$

\noindent
Now the $\mbox{SO}(4p)$--matrices which commute with each of the complex structures in the hypercomplex structure $\mQ$ on $\mR^{4p}$, form a subgroup of $\mbox{SO}_{\mI}(4p)$, denoted by $\mbox{SO}_{\mQ}(4p)$, which is isomorphic with the symplectic group $\mbox{Sp}( p)$.
Recall that the symplectic group $\mbox{Sp}(p)$ is the real Lie group of quaternion $p \times p$ matrices preserving the symplectic inner product $\langle \xi , \eta \rangle_{\mH} = \xi_1 \overline{\eta}_1 + \xi_2 \overline{\eta}_2 + \cdots +  \xi_p \overline{\eta}_p$, $\xi, \eta \in \mH^p$, or equivalently,
$$
\mbox{Sp}(p) = \left \{ A \in \mbox{GL}_p(\mH) : AA^\ast = E_p \right \}
$$  
Quite naturally, the subgroup $\mbox{SO}_{\mQ}(4p)$ of $\mbox{SO}(4p)$ has a double covering by $\mbox{Spin}_{\mQ}(4p)$, the subgroup of $\mbox{Spin}(4p)$ consisting of the $\mbox{Spin}(4p)$--elements which commute with both $s_{\mI}$ and $s_{\mJ}$, where now $s_{\mJ}$ is the $\mbox{Spin}(4p)$--element corresponding to the complex structure $\mJ_p$. Recall, see  (\ref{sI}), that $s_\mI$, corresponding to $\mI_{4p}$, is given by $s_{\mI} = s_1 \cdots s_{2p}$, where $s_j = \frac{\sqrt{2}}{2} \bigl ( 1 + e_{2j-1} e_{2j} \bigr )$, $j=1,\ldots, 2p$. Similarly, for $s_{\mJ}$ we find
\begin{equation}
s_{\mJ} = \widetilde{s_1} \cdots \widetilde{s_p}, \quad \widetilde{s_j} = \frac{1}{2} \bigl ( 1 + e_{4j-3} e_{4j-1} \bigr ) \bigl (1 - e_{4j-2} e_{4j} \bigr ), \; j=1,\ldots, p
\label{sJ}
\end{equation}
For the corresponding picture at the level of the Lie algebras we refer to \cite{paper1}.\\[-2mm]

The introduction of a hypercomplex structure leads to a function theory for so-called {\em quaternionic Clifford analysis}, where the fundamental invariance will be that of the symplectic group $\mbox{Sp}(p)$. The most genuine way to introduce the new concept of quaternionic monogenicity is to directly generalize the system (\ref{hmon}), expressing hermitian monogenicity, now making use of the additional rotated Dirac operators $\dirac_\mJ = \mJ_{4p}[\dirac]$ and $\dirac_\mK = \mK_{4p}[\dirac]$, whence the following definition.

\begin{definition}
\label{defqmon}
A differentiable function $F: \mR^{4p} \longrightarrow \mS$ is called quaternionic monogenic (or q--monogenic for short) in some region $\Omega$ of $\mR^{4p}$, if and only if in that region $F$ is a solution of the system
\begin{equation}
\dirac F = 0, \quad \dirac_\mI F = 0, \quad \dirac_\mJ F = 0, \quad \dirac_\mK F = 0
\label{qmon}
\end{equation}
\end{definition}

\noindent
Observe that, in a similar way as it was possible to introduce hermitian monogenicity without involving complex numbers, the above definition expresses q--monogenicity without having to resort to quaternions.\\[-2mm]

Quite naturally, there is an alternative characterization of q--monogenicity in terms of the hermitian Dirac operators, yet still not involving quaternions. We recall these hermitian Dirac operators in the actual dimension:
\begin{eqnarray*}
\upz & = & \sum_{k=1}^{2p} \p_{z_k} \gfd_k \ = \ \sum_{j=1}^p ( \p_{z_{2j-1}} \gfd_{2j-1} + \p_{z_{2j}} \gfd_{2j} ) \\
\upzd & = & \sum_{k=1}^{2p} \p_{\overline{z}_k} \gf_k \ = \ \sum_{j=1}^p ( \p_{\overline{z}_{2j-1}} \gf_{2j-1} + \p_{\overline{z}_{2j}} \gf_{2j} )
\end{eqnarray*}
and determine their images under the action of $\mJ_{4p}$:
\begin{eqnarray*}
\upzJ \ = \ \mJ_{4p} [\upz] & = &  \sum_{j=1}^p ( \p_{z_{2j}} \gf_{2j-1} - \p_{z_{2j-1}} \gf_{2j} ) \\
\upzJd \ = \ \mJ_{4p} [\upzd] & = & \sum_{j=1}^p ( \p_{\overline{z}_{2j}} \gfd_{2j-1} - \p_{\overline{z}_{2j-1}} \gfd_{2j} )
\end{eqnarray*}
Similarly we can introduce the variables
\begin{eqnarray*}
\uzJ \ = \ \mJ_{4p} [\uz] & = & \sum_{j=1}^p ( z_{2j} \, \gfd_{2j-1} - z_{2j-1} \, \gfd_{2j} ) \\
\uzJd \ = \ \mJ_{4p} [\uzd] & = & \sum_{j=1}^p ( \overline{z}_{2j} \, \gf_{2j-1} - \overline{z}_{2j-1} \, \gf_{2j} )
\end{eqnarray*}
Here the formulae $\mJ_{4p} [ \gf_{2j-1} ] = - \gfd_{2j}$, $\mJ_{4p} [\gf_{2j}] = \gfd_{2j-1}$, $\mJ_{4p}[\gfd_{2j-1}] = - \gf_{2j}$ and $\mJ_{4p}[\gfd_{2j}] = \gf_{2j-1}$ were used.\\[-2mm]

Now the Dirac operator $\dirac$ and its rotated versions $\dirac_\mI$, $\dirac_\mJ$, $\dirac_\mK$ may be expressed in terms of the hermitian Dirac operators $( \upz,\upzd)$ and their rotated versions $(\upzJ,\upzJd)$. We indeed have
\begin{align*}
\upz &= \pi^+[\dirac] = \frac{1}{4} ( {\bf 1} + i \, \mI_{4p} ) [\dirac] = \frac{1}{4} \left( \dirac + i \, \dirac_\mI \right)\\
\upzd &= \pi^-[\dirac] = -\frac{1}{4} ( {\bf 1} - i \, \mI_{4p} ) [\dirac] = \frac{1}{4} \left( - \dirac + i \, \dirac_\mI\right)\\
\upzJ &= \mJ_{4p} [\pi^+[\dirac]] = \frac{1}{4} ( \mJ_{4p} + i \, \mK_{4p}) [\dirac] = \frac{1}{4} \left( \dirac_\mJ + i \, \dirac_\mK \right)\\
\upzJd &=  \mJ_{4p} [\pi^-[\dirac]] = -\frac{1}{4} ( \mJ_{4p} - i \, \mK_{4p}) [\dirac] = \frac{1}{4} \left( \dirac_\mJ - i \, \dirac_\mK \right)
\end{align*}
whence, conversely,
$$
\dirac = 2 ( \upz - \upzd), \quad i \, \dirac_\mI = 2 ( \upz + \upzd ), \quad \dirac_\mJ = 2 ( \upzJ - \upzJd ), \quad  i \, \dirac_\mK =2 ( \upzJ + \upzJd) 
$$
This leads to an alternative characterization of q--monogenicity in the following proposition.
\begin{proposition}
A differentiable function $F: \mR^{4p} \simeq \mC^{2p} \longrightarrow \mS$ is q--monogenic in the region $\Omega \subset \mR^{4p}$ if and only if $F$ is in $\Omega$ a simultaneous null solution of the operators $\upz$, $\upzd$, $\upzJ$ and $\upzJd$.
\end{proposition}

As the identification of the underlying symmetry group is necessary for the further development of the function theory, the next result is crucial.
\begin{proposition}
The operators $\upz$, $\upzd$, $\upzJ$ and $\upzJd$ are invariant under the action of the symplectic group $\mbox{\em Sp}(p)$.
\end{proposition}

\pf
The action of a $\mbox{Spin}(4p)$--element $s$ on a spinor--valued function $F$ is the so-called $L$--action given by $L(s) [ F(\uX)] = s F ( s^{-1} \uX s)$. The Dirac operator $\dirac$ is invariant under $\mbox{Spin}(4p)$, i.e.\ $[L(s),\dirac] = 0$, for all $s \in \mbox{Spin}(4p)$, which can be explained by
$$
L(s) \dirac_{\uX} F(\uX) = s \dirac_{s^{-1} \uX s} F (s^{-1} \uX s) = s ( s^{-1} \dirac_{\uX} s ) F ( s^{-1} \uX s)  = \dirac_{\uX} L(s) F(\uX)
$$
Recall that $\mbox{Sp}(p)$ is isomorphic to a subgroup of $\mbox{Spin}(4p)$, whence the Dirac operator $\dirac$ is, quite trivially, invariant under the action of $\mbox{Sp}(p)$. The invariance of the operators $\upz$, $\upzd$, $\upzJ$ and $\upzJd$ now follows from the fact that their respective definitions involve projection operators which are commuting with the $\mbox{Sp}(p)$--elements.
\qed


\section{A further decomposition of spinor space}
\label{spinordecomposition}


The main aim of this paper being the decomposition of the space $\mcP(\mR^{4p};\mS)$ of spinor--valued polynomials into Sp$( p)$--irreducibles, we should first take care of the irreducibility of the value space as an Sp$( p)$--representation. Spinor space $\mS$, which was already decomposed into U$(2p)$--irreducible homogeneous parts $\mS^r$, should further be decomposed into Sp$( p)$--irreducibles, which we will call {\em symplectic cells}. We sketch this decomposition and refer to \cite{paper1} for the details.\\[-2mm]

\noindent
First we introduce the Sp$( p)$--invariant left multiplication operators
$$
P = \gf_2 \gf_1 + \gf_4 \gf_3 + \ldots + \gf_{2p} \gf_{2p-1}, \qquad
Q = \gfd_1 \gfd_2 + \gfd_3 \gfd_4 + \ldots + \gfd_{2p-1} \gfd_{2p} \ = \ P^{\dagger}
$$
for which $P: \mS^r \rightarrow \mS^{r-2}$ and $Q: \mS^r \rightarrow \mS^{r+2}$. 
Together with the spin--Euler operator $\beta$, they generate an $\gsl_2(\mC)$--structure as is seen from the following relations which may be directly verified. 

\begin{lemma}
\label{lemmaS1}
One has
\begin{enumerate}
\item[(i)] $[P,Q] = p - \beta$;
\item[(ii)] $[p-\beta,P] = 2P$;
\item[(iii)] $[p-\beta,Q] = -2Q$
\end{enumerate}
\end{lemma}

\noindent Also the following results are straightforward.
\begin{lemma}
\label{lemmaS3}
One has
\begin{enumerate}
\item[(i)] $\mbox{\em Ker} \, P|_{\mS^r} = \{0\}$ for $r=p+1,\ldots,2p$;
\item[(ii)] $\mbox{\em Ker} \, Q|_{\mS^r} = \{0\}$ for $r=0,\ldots,p-1$;
\item[(iii)] $\mbox{\em Ker} \, P|_{\mS^p} = \mbox{\em Ker} \, Q|_{\mS^p}$. 
\end{enumerate}
\end{lemma}

Next we define, for $r=0,\ldots,p$, the subspaces
$$
\mS_r^r = \mbox{Ker} \, P |_{\mS^r}, \qquad \mS_r^{2p-r} = \mbox{Ker} \, Q |_{\mS^{2p-r}}
$$
and for $k=0,\ldots,p-r$, the subspaces
$$
\mS_r^{r+2k} = Q^k \, \mS^r_r, \qquad \mS_r^{2p-r-2k} = P^k \, \mS^{2p-r}_r
$$

\begin{lemma}
\label{lemmaS4}
One has
$$
\mbox{\em Ker} \, P|_{\mS} = \bigoplus_{r=0}^p \mS_r^r, \qquad
\mbox{\em Ker} \, Q|_{\mS} = \bigoplus_{r=0}^p \mS_r^{2p-r} 
$$
and, for $k= 0,1,\ldots,p-r-1$,
\begin{itemize}
\item $Q$ is an isomorphism $\mS_r^{r+2k} \longrightarrow \mS_r^{r+2k+2}$ with inverse $Q^{-1} = \displaystyle\frac{1}{\alpha_r^k} \, P$;\\[-4mm]
\item  $P$ is an isomorphism $\mS_r^{2p-r-2k} \longrightarrow \mS_r^{2p-r-2k-2}$ with inverse $P^{-1} = \displaystyle\frac{1}{\alpha_r^{p-r-k-1}} \, Q$
\end{itemize}
where the coefficients $\alpha_r^k$ are given by
$$
\alpha_r^k = (k+1)(p-r-k) = \alpha_r^{p-r-k-1}
$$
which implies that the composition of the multiplicative operators $P$ and $Q$ is constant on each symplectic cell; more specifically one has
\begin{itemize}
\item $P \, Q = \alpha_r^k$ on \, $\mS_r^{r+2k}$ and on \, $\mS_r^{2p-r-2k-2}$
\item $Q \, P = \alpha_r^k$ on \,  $\mS_r^{r+2k+2}$ and on \, $\mS_r^{2p-r-2k}$
\end{itemize}
\end{lemma}

\begin{proposition}
\label{propositionS1}
One has, for all $r=0,\ldots,p$,
$$
\mS^r =\bigoplus_{j=0}^{\lfloor \frac{r}{2} \rfloor} \mS_{r-2j}^r, \qquad
\mS^{2p-r} =\bigoplus_{j=0}^{\lfloor \frac{r}{2} \rfloor} \mS_{r-2j}^{2p-r}
$$
and each of the symplectic cells $\mS_s^r$ in the above decompositions is an irreducible $\mbox{Sp}(p)$--representation. 
\end{proposition}

\noindent
The above decomposition of the homogeneous spinor subspaces into symplectic cells is depicted in the triangular scheme below. The vertical columns correspond to the homogeneous spinor subspaces; they are irreducible representations of U$(2p)$. Each of the symplectic cells on one particular row is a realization of one and the same irreducible Sp$( p)$ representation. The operators $P$ and $Q$ allow for  shifting horizontally between those isomorphic cells, and, as was pointed out above (see Lemma \ref{lemmaS4}), on each symplectic cell the products $PQ$ and $QP$ act as constants.\\

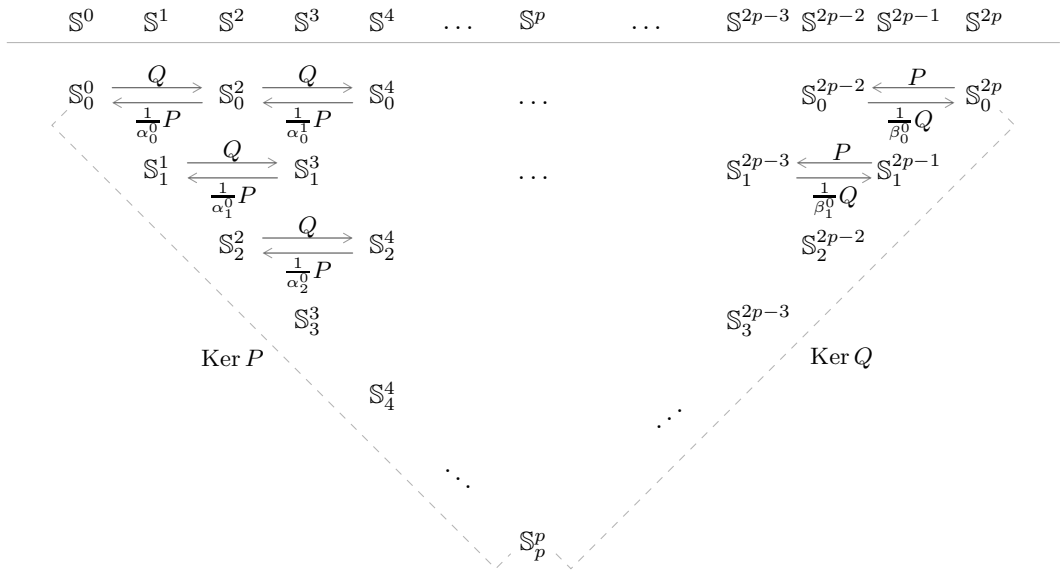
\begin{center}
\begin{tikzpicture}[scale=1]
\node at (0,0) {$\mS^0$};
\node at (1,0) {$\mS^1$};
\node at (2,0) {$\mS^2$};
\node at (3,0) {$\mS^3$};
\node at (4,0) {$\mS^4$};
\node at (5,-0.1) {$\ldots$};
\node at (6,0) {$\mS^p$};
\node at (7.5,-0.1) {$\ldots$};
\node at (9,0) {$\mS^{2p-3}$};
\node at (10,0) {$\mS^{2p-2}$};
\node at (11,0)  {$\mS^{2p-1}$};
\node at (12,0) {$\mS^{2p}$};
\draw[-, draw opacity=0.3] (-1,-0.3) -- (13,-0.3);
\node at (0,-1) {$\mS^0_0$};
\node at (2,-1) {$\mS^2_0$};
\node at (4,-1) {$\mS^4_0$};
\node at (6,-1.1) {$\ldots$};
\node at (6,-2.1) {$\ldots$};
\node at (10,-1) {$\mS^{2p-2}_0$};
\node at (12,-1) {$\mS^{2p}_0$};
\node at (1,-2) {$\mS^1_1$};
\node at (3,-2) {$\mS^3_1$};
\node at (11,-2) {$\mS^{2p-1}_1$};
\node at (9,-2) {$\mS^{2p-3}_1$};
\node at (2,-3) {$\mS_2^2$};
\node at (4,-3) {$\mS_2^4$};
\node at (10,-3) {$\mS_2^{2p-2}$};
\node at (9,-4) {$\mS_3^{2p-3}$};
\node at (3,-4) {$\mS_3^3$};
\node at (4,-5) {$\mS_4^4$};
\node at (5,-6) {$\ddots$};
\node at (6,-7) {$\mS_p^p$};
\node at (7.8,-5.2) {$\iddots$};
\draw[->,>=angle 60,draw opacity=0.4] (0.4,-0.9) -- (1.6,-0.9);
\node at (1,-0.75) {\small $Q$};
\draw[<-,>=angle 60,draw opacity=0.4] (0.4,-1.1) -- (1.6,-1.1);
\draw[<-,>=angle 60,draw opacity=0.4] (10.5,-0.9) -- (11.6,-0.9);
\draw[->,>=angle 60,draw opacity=0.4] (10.45,-1.1) -- (11.6,-1.1);
\node at (11.1,-0.75) {\small $P$};
\node at (3,-0.75) {\small $Q$};
\draw[->,>=angle 60,draw opacity=0.4] (2.4,-0.9) -- (3.6,-0.9);
\draw[<-,>=angle 60,draw opacity=0.4] (2.4,-1.1) -- (3.6,-1.1);
\node at (2,-1.75) {\small $Q$};
\draw[->,>=angle 60,draw opacity=0.4] (1.4,-1.9) -- (2.6,-1.9);
\draw[<-,>=angle 60,draw opacity=0.4] (1.4,-2.1) -- (2.6,-2.1);
\draw[<-,>=angle 60,draw opacity=0.4] (9.5,-1.9) -- (10.5,-1.9);
\draw[->,>=angle 60,draw opacity=0.4] (9.5,-2.1) -- (10.5,-2.1);
\node at (10.1,-1.75) {\small $P$};
\node at (3,-2.75) {\small $Q$};
\draw[->,>=angle 60,draw opacity=0.4] (2.4,-2.9) -- (3.6,-2.9);
\draw[<-,>=angle 60,draw opacity=0.4] (2.4,-3.1) -- (3.6,-3.1);
\draw[draw opacity = 0.3,dashed] (-0.2,-1.2) -- (-0.4,-1.4) -- (5.5,-7.3) -- (5.7,-7.1) ;
\node at (2,-4.5) {\small $\mbox{Ker} \, P$};
\draw[draw opacity = 0.3,dashed] (6.3,-7.1) -- (6.5,-7.3) -- (12.4,-1.4) -- (12.2,-1.2);
\node at (10.1,-4.5) {\small $\mbox{Ker} \, Q$};
\node at (1,-1.4) {\small $\frac{1}{\alpha_0^0} P$};
\node at (3,-1.4) {\small $\frac{1}{\alpha_0^1} P$};
\node at (2,-2.4) {\small $\frac{1}{\alpha_1^0} P$};
\node at (3,-3.4) {\small $\frac{1}{\alpha_2^0} P$};
\node at (11,-1.4) {\small $\frac{1}{\beta_0^0} Q$};
\node at (10,-2.4) {\small $\frac{1}{\beta_1^0} Q$};
\end{tikzpicture}
\end{center}

\noindent


\section{$\gosp(4|2)$--monogenic functions}
\label{osp42monogenic}


According to the splitting of spinor space into symplectic cells,  a function $F : \mR^{4p} \longrightarrow \mS$ can be decomposed into components taking values in these symplectic cells:
$$
F = \sum_{r=0}^n \, F^r =  \sum_{r=0}^n \, \sum_s F^r_s, \qquad F^r_s : \mR^{4p} \longrightarrow \mS^r_s
$$
Expressing that the function $F$ is q--monogenic is, quite surprisingly, equivalent with the q--monogenicity of each of the components $F^r_s$, and leads to systems of first order differential equations the form of which  depends on where the symplectic cell is situated in the triangular decomposition scheme of spinor space. For a detailed study of those systems of equations we refer to \cite{paper2}. One could expect now the building blocks of the Fischer decomposition under the Sp$( p)$--action to be the spaces $\mcQ^{r,s}_{a,b}$ of quaternionic monogenic bi--homogeneous polynomials with values in a symplectic cell:
$$
\mcQ^{r,s}_{a,b} = \mcP^{r,s}_{a,b} \cap \mathrm{Ker}  \left(\upz, \upzd, \upzJ, \upzJd \right) = \mcP_{a,b}(\mR^{4p};\mS^r_s) \cap \mathrm{Ker}  \left(\upz, \upzd, \upzJ, \upzJd \right)
$$ 
Unfortunately these spaces  $\mcQ^{r,s}_{a,b}$ are reducible under the action of the group Sp$( p)$, whence the notion of q--monogenicity has to be refined. In order to determine the suitable refinement, the traditional path:\\[-2mm]

\noindent ''differential operator(s)" $\rightarrow$ "invariance group"  $\rightarrow$ "Fischer decomposition" $\rightarrow$ "Howe dual pair"\\[-2mm]

\noindent has to be abandoned. Instead, we will first turn our attention to the Howe dual partner of the chosen invariance group $G$. Indeed, if the Lie (super)algebra $\geg$ is the Howe dual partner of $G$, to be found within the Weyl algebra of polynomial differential operators, the characterizing differential operators, which are $G$--invariant, correspond to the {\em negative roots} in the root system obtained via the $ad(\gH)$--action of the Cartan subalgebra $\gH$. In the present case where $G = \mathrm{Sp}( p)$, we want the already established differential operators $\upz$, $\upzd$, $\upzJd$ and $\upzJd$ and their algebraic counterparts $\uz$, $\uzd$, $\uzJd$ and $\uzJd$, which  indeed all are Sp$( p)$--invariant, to belong to (the odd part of) the Lie (super)algebra $\geg$, which has to be closed under the Lie super bracket. 
Computing the anti--commutators of those differential and multiplication operators (see the complete list below) we find, next to the expected $\gsl(2,\mC)$ generators, viz. $\mE_{\uz} + \mE_{\uzd} + 2p$, $r^2$, $\Delta_{4p}$, and the operator $\mE_{\uz} - \mE_{\uzd}$, which already appeared in the context of hermitian Clifford analysis, and next to the shifting operators $P$ and $Q$ used in the definition of the symplectic cells, the new scalar differential operators
$$
\SE =  \sum_{k=1}^p  \, z_{2k-1} \, \p_{\olz_{2k}} - z_{2k} \, \p_{\olz_{2k-1}} \quad \mathrm{and} \quad \SED  = \sum_{k=1}^p  \, \olz_{2k-1} \, \p_{z_{2k}} - \olz_{2k} \, \p_{z_{2k-1}}
$$
enjoying the following properties (see \cite{paper3}).

\begin{lemma}
The operators $\SE$ and $\SED$ are invariant under the symplectic action.
\end{lemma}

\begin{lemma}
One has
$$
\gsl(2,\mC) \cong {\rm Alg}_\mC \left( \mE_{\uz} - \mE_{\uz}^\dagger , \SE, \SED  \right)
$$
\end{lemma}

\begin{lemma}
The three $\gsl(2,\mC)$--structures  $\mathrm{Alg}_{\mC} \left( \mE_{\uz} + \mE_{\uz}^\dagger + 2p, \onehalf r^2, - \onehalf \Delta_{4p} \right)$, $\mathrm{Alg}_{\mC} \left( p-\beta, P, Q  \right)$ and ${\mathrm Alg}_\mC \left( \mE_{\uz} - \mE_{\uz}^\dagger , \SE, \SED  \right)$ are mutually commuting.

\end{lemma}

So it becomes clear that the Howe dual partner to the Lie group Sp$( p)$ is the Lie superalgebra
$$
\mathfrak{osp}(4|2) = \geg_0 \oplus \geg_1 = \left( \gso(4) \oplus \gsp(2) \right) \oplus \geg_1 = \left( \gsl(2) \oplus \gsl(2) \oplus \gsl(2) \right) \oplus \geg_1
$$
with
$$
\gsl(2) \cong \mathrm{Alg}_{\mC} \left( \mE_{\uz} + \mE_{\uz}^\dagger + 2p, \onehalf \nz^2, - \onehalf \Delta_{4p} \right)
\cong \mathrm{Alg}_{\mC} \left( \mE_{\uz} - \mE_{\uz}^\dagger, \SE, \SED \right)
\cong \mathrm{Alg}_{\mC} \left( p-\beta, P, Q  \right)
$$
and
$$
\geg_1 \cong \mathrm{span}_{\mC} \left( \uz, \uzd, \upz, \upzd  \right) \oplus \mathrm{span}_{\mC} \left( \uzJ, \uzJd, \upzJ, \upzJd  \right)
$$
The $\mZ$--gradation of the $17$ dimensional orthosymplectic Lie superalgebra $\mathfrak{osp}(4|2)$ has the structure
$$
\mathfrak{osp}(4|2) = \mcG_{-2} \oplus  \mcG_{-1} \oplus \mcG_{0} \oplus \mcG_{1} \oplus \mcG_{2}
$$
where
$$
\mcG_{-2} = \mathrm{Alg}_{\mC} \left( \nz^2 \right)
$$
$$
\mcG_{-1} = \mathrm{span}_{\mC} \left( \uz, \uzd, \uzJ, \uzJd \right)
$$
$$
\mcG_{0} = \mathfrak{H}  \cup \{ \SE, \SED, P, Q \}
$$
$$
\mcG_{1} = \mathrm{span}_{\mC} \left( \upz, \upzd, \upzJ, \upzJd  \right)
$$
$$
\mcG_{2} = \mathrm{Alg}_{\mC} \left( \Delta_{4p} \right)
$$
with 
$$\mathfrak{H} = \mathrm{Alg}_{\mC} \left( \mE_{\uz} + \mE_{\uz}^\dagger + 2p,  \mE_{\uz} - \mE_{\uz}^\dagger, p-\beta \right)
$$
the Cartan subalgebra. Note that the $8$ dimensional odd part $\geg_1$ of $\mathfrak{osp}(4|2)$ is nothing but the spinorial representation of the even subalgebra $\geg_0$, i.e. $\mC^2 \otimes \mC^2 \otimes \mC^2$, such that each odd element has a label in $\mZ^3$ consisting of the eigenvalues under the action of the Cartan subalgebra, i.e.
$$
\uz \rightarrow (1,1,1), \quad \uzd \rightarrow (1,-1,-1), \quad \upz \rightarrow (-1,-1,-1) , \quad \upzd \rightarrow (-1,1,1)
$$
$$
\uzJ \rightarrow (1,1,-1), \quad \uzJd \rightarrow (1,-1,1), \quad \upzJ \rightarrow (-1,-1,1), \quad \upzJd \rightarrow (-1,1,-1)
$$
as also shows from the following overview of the (non--trivial) (anti--)commutation relations:\\

\noindent $\bullet$ within $\mcG_0$
{\scriptsize
\begin{align*}
[\mE_{\uz} + \mE_{\uz}^\dagger + 2p, \SE ] &=  0 & [\mE_{\uz} + \mE_{\uz}^\dagger + 2p, \SED ] &= 0\\
[  \mE_{\uz} - \mE_{\uz}^\dagger, \SE] &= 2 \SE & [  \mE_{\uz} - \mE_{\uz}^\dagger, \SED] &= - 2 \SED\\
[ p - \beta, P ] &= 2P & [ p - \beta, Q ] &= - 2Q\\
[ \SE , \SED ] &=  \mE_{\uz} - \mE_{\uz}^\dagger & [ P , Q ] &= p - \beta
\end{align*}
}
\vspace*{-2mm}

$\bullet$ between $\mcG_0$ and $\mcG_1$
{\scriptsize
\begin{align*}
[\mE_{\uz} + \mE_{\uz}^\dagger + 2p, \upz ] &=  -\upz  & [\mE_{\uz} + \mE_{\uz}^\dagger + 2p, \upzd ] &=  -\upzd &[\mE_{\uz} + \mE_{\uz}^\dagger + 2p, \upzJ ] &= -\upzJ  &[\mE_{\uz} + \mE_{\uz}^\dagger + 2p, \upzJd ] &= -\upzJd\\
[\mE_{\uz} - \mE_{\uz}^\dagger, \upz ] &=  -\upz  & [\mE_{\uz} - \mE_{\uz}^\dagger, \upzd ] &= \upzd  & [\mE_{\uz} - \mE_{\uz}^\dagger, \upzJ ] &= -\upzJ  & [\mE_{\uz} - \mE_{\uz}^\dagger, \upzJd ] &= \upzJd\\
[p-\beta, \upz ] &= - \upz  & [p-\beta, \upzd ] &=  \upzd & [p-\beta, \upzJ ] &=  \upzJ &[p-\beta, \upzJd ] &= - \upzJd\\
[\SE, \upz ] &=  - \upzJd & [\SE, \upzd ] &= 0  & [\SE, \upzJ ] &=  \upzd & [\SE, \upzJd ] &= 0\\
[\SED, \upz ] &=  0 & [ \SED, \upzd ] &= \upzJ &[\SED, \upzJ ] &= 0 &[\SED, \upzJd ] &= - \upz\\
[P, \upz ] &=  -\upzJ &[P, \upzd ] &=  0 &[P, \upzJ ] &= 0 &[P, \upzJd ] &= \upzd\\
[Q, \upz ] &=   0 &[Q, \upzd ] &=  \upzJd & [Q, \upzJ ] &=  - \upz &[Q, \upzJd ] &= 0
\end{align*}
}
\vspace*{-2mm}

$\bullet$ between $\mcG_0$ and $\mcG_{-1}$
{\scriptsize 
\begin{align*}
[\mE_{\uz} + \mE_{\uz}^\dagger + 2p, \uz ] &=  \uz  & [\mE_{\uz} + \mE_{\uz}^\dagger + 2p, \uzd ] &=  \uzd & [\mE_{\uz} + \mE_{\uz}^\dagger + 2p, \uzJ ] &= \uzJ  & [\mE_{\uz} + \mE_{\uz}^\dagger + 2p, \uzJd ] &= \uzJd\\
[\mE_{\uz} - \mE_{\uz}^\dagger, \uz ] &= \uz  & [\mE_{\uz} - \mE_{\uz}^\dagger, \uzd ] &=  -\uzd & [ \mE_{\uz} - \mE_{\uz}^\dagger, \uzJ ] &=  \uzJ & [\mE_{\uz} - \mE_{\uz}^\dagger, \uzJd ] &= -\uzJd\\
[p-\beta, \uz ] &=   \uz & [p-\beta, \uzd ] &= - \uzd  & [p-\beta, \uzJ ] &=  - \uzJ & [p-\beta, \uzJd ] &= \uzJd\\
[\SE, \uz ] &=  0 & [\SE, \uzd ] &= -\uzJ & [\SE, \uzJ ] &=  0 & [\SE, \uzJd ] &= \uz\\
[\SED, \uz ] &=  \uzJd & [\SED, \uzd ] &= 0 & [\SED, \uzJ ] &= -\uzd & [\SED, \uzJd ] &= 0\\
[P, \uz ] &=  0 & [P, \uzd ] &=  -\uzJd & [P, \uzJ ] &= \uz & [P, \uzJd ] &= 0\\
[Q, \uz ] &=   \uzJ & [Q, \uzd ] &=  0 & [Q, \uzJ ] &=  0 & [Q, \uzJd ] &= -\uzd
\end{align*}
}
\vspace*{-2mm}

$\bullet$ between $\mcG_0$ and $\mcG_{2}$
{\scriptsize \begin{align*}
[\mE_{\uz} + \mE_{\uz}^\dagger + 2p, \Delta_{4p} ] &=  -2 \Delta_{4p} & [\SE, \Delta_{4p} ] &=  0 \\
[\mE_{\uz} - \mE_{\uz}^\dagger,\Delta_{4p} ] &= 0 & [\SED,\Delta_{4p} ] &= 0
\end{align*}
}
\vspace*{-2mm}

$\bullet$ between $\mcG_0$ and $\mcG_{-2}$
{\scriptsize
\begin{align*}
[\mE_{\uz} + \mE_{\uz}^\dagger + 2p, \nz^2 ] &=  2 \nz^2 & [\SE, \nz^2 ] &= 0\\
[\mE_{\uz} - \mE_{\uz}^\dagger,  \nz^2 ] &= 0  & [\SED, \nz^2 ] &= 0
\end{align*}
}
\vspace*{-2mm}

$\bullet$ within $\mcG_{1}$
{\scriptsize
\begin{align*}
 \{ \upz, \upzd \} &=  \onequarter \Delta_{4p}  & \{\upz, \upzJ \} &= 0  & \{ \upz, \upzJd \} = 0\\
\{\upzd, \upzJ \} &= 0  & \{\upzd, \upzJd \} &= 0 & \{\upzJ, \upzJd \} &= \onequarter \Delta_{4p}
\end{align*}
}
\vspace*{-2mm}

$\bullet$ between $\mcG_{1}$ and $\mcG_{-1}$ 
{\scriptsize
\begin{align*}
\{\upz, \uz \} &=  \mE_{\uz}  + \beta & \{\upz, \uzd \} &= 0  & \{\upz, \uzJ\} &= -2Q & \{\upz, \uzJd \} &= \SED\\
\{\upzd, \uz\} &=  0 & \{\upzd, \uzd \} &= \mE_{\uz}^\dagger + 2p - \beta & \{\upzd, \uzJ \} &=  - \SE & \{\upzd, \uzJd \} &= 2P\\
\{\upzJ, \uz\} &=  -2P & \{\upzJ, \uzd \} &= - \SED & \{\upzJ, \uzJ \} &= \mE_{\uz} + 2p - \beta & \{\upzJ, \uzJd\} &= 0\\
\{\upzJd, \uz \} &=  \SE & \{\upzJd, \uzd \} &=  2Q & \{\upzJd, \uzJ \} &= 0 & \{\upzJd, \uzJd \} &= \mE_{\uz}^\dagger + \beta
\end{align*}
}
\vspace*{-2mm}

$\bullet$ the operators in $\mcG_{1}$ and $\mcG_{2}$ are commuting

$\bullet$ between $\mcG_{1}$ and $\mcG_{-2}$ 
{\scriptsize
\begin{align*}
[\upz, \nz^2 ] &=  \uzd & [\upzJ, \nz^2 ] &= \uzJd\\
[\upzd, \nz^2 ] &= \uz & [\upzJd, \nz^2 ] &= \uzJ
\end{align*}
}
\vspace*{-2mm}

$\bullet$ within $\mcG_{-1}$
{\scriptsize
\begin{align*}
 \{\uz, \uzd \} &= \nz^2  & \{\uz, \uzJ\} &= 0 & \{\uz, \uzJd \} &= 0\\
\{\uzd, \uzJ \} &=  0 & \{\uzd, \uzJd \} &= 0 &  \{\uzJ, \uzJd \} &= \nz^2
\end{align*}
}
\vspace*{-2mm}

$\bullet$ between $\mcG_{-1}$ and $\mcG_{2}$ 
{\scriptsize
\begin{align*}
[\uz, \Delta_{4p} ] &=  - 4 \upzd & [\uzJ, \Delta_{4p} ] &=  - 4 \upzJd\\
[\uzd,\Delta_{4p} ] &=  - 4 \upz & [\uzJd,\Delta_{4p} ] &= - 4 \upzJ
\end{align*}
}
\vspace*{-2mm}

$\bullet$ the operators in $\mcG_{-1}$ and $\mcG_{-2}$ are commuting

$\bullet$ between $\mcG_{2}$ and $\mcG_{-2}$ 
{\scriptsize
\begin{align*}
[ \Delta_{4p} , \nz^2 ] &= 4 \left( \mE_{\uz} + \mE_{\uz}^\dagger + 2p  \right)
\end{align*}
}

\noindent In view of the above commutation relations between $\mcG_0$ and $\mcG_1$ we have the following result.

\begin{proposition}
If the function $F$ is q--monogenic, then also the functions $\SE F$, $\SED F$, $PF$ and $QF$ are q--monogenic.
\end{proposition}

The Howe dual pair $\gosp(4|2)  \times \mathrm{Sp}( p)$ then inevitably leads to the following necessary refinement of the concept of quaternionic monogenicity.

\begin{definition} 
\label{defsqmon}
A differentiable function $F : \mR^{4p} \longrightarrow \mS$ is called $\gosp(4|2)$--monogenic in some region $\Omega$ of $\mR^{4p}$, if and only if in that region $F$ is q--monogenic and at the same time a null--solution of both the operators $\mcE$ and $P$.
\end{definition}

\noindent
Quite naturally, there are similar definitions involving the operators $\SED$ and $Q$. It is clear that a function which is in Ker $P$ (respectively Ker $Q$) only can take its values in the symplectic cells situated on the left (respectively right) edge of the triangular spinor scheme; the condition to belong to either Ker $\SE$ or Ker $\SED$ has implications on the bidegree of the polynomials considered. This is taken into account in the definition of the following polynomial function spaces, which will appear in our Fischer decomposition aimed at.

\begin{definition}
If $a \geq b$ then 
\begin{align*}
\mcS^{r}_{a,b} &= \mcQ^{r,r}_{a,b} \cap \mathrm{Ker} \, (\SE,P), \ \ r=0,\ldots,p\\
\mcT^{r}_{a,b} &= \mcQ^{2p-r,r}_{a,b} \cap \mathrm{Ker} \, (\SE,Q), \ \ r=p,\ldots,2p.
\end{align*}
If $a \leq b$ then 
\begin{align*}
\mcS^{r \dagger}_{a,b} &= \mcQ^{r,r}_{a,b} \cap \mathrm{Ker} \, (\SED,P), \ \ r=0,\ldots,p\\
\mcT^{r \dagger}_{a,b} &= \mcQ^{2p-r,r}_{a,b} \cap \mathrm{Ker} \, (\SED,Q), \ \  r=p,\ldots,2p.
\end{align*}
\end{definition}

\begin{remark}
As in \cite{paper3} it may be shown that if $a > b$ then $\mcQ^{r,r}_{a,b} \cap \mathrm{Ker} \, \SED = \{0\}$, while if $a < b$ then  $\mcQ^{r,r}_{a,b} \cap \mathrm{Ker} \, \SE = \{0\}$. Moreover $\mcS^{r}_{a,a} = \mcS^{r \dagger}_{a,a}$ and   $\mcT^{r}_{a,a} = \mcT^{r \dagger}_{a,a}$.
\end{remark}


\section{Fischer decomposition into Sp$( p)$--irreducibles}
\label{fischer}


In this section we will eventually prove (see Theorem \ref{finaltheo}) a Fischer decomposition of the space $\mcP(\mR^{4p};\mS)$ of all spinor--valued polynomials in terms of spaces of spherical harmonics, each of which is the product of a space $\mcS^{r}_{a,b}$  of spherical $\gosp(4|2)$--monogenics, introduced in the preceding section, with an appropriate embedding factor involving the variables $\uz, \uzd, \uzJ, \uzJd$ and the operators $\SED$ and $Q$. We will prove this result through a Fischer decomposition of spaces of spherical harmonics, since spaces of polynomials may always be decomposed into spaces of homogeneous harmonic  polynomials with values in appropriate invariant subspaces of spinor space.  More precisely we will first concentrate on obtaining a Fischer decomposition of the space $\mcH_{a,b}(\mR^{4p};\mS)$ of spinor--valued bi--homogeneous harmonic  polynomials in terms of  $\mathrm{Sp}( p)$--irreducibles. \\[-2mm]

A first important step was already taken in \cite{paper3}, where we obtained a Fischer decomposition of the space $\mcH_{a,b}(\mR^{4p};\mC)$ of scalar--valued bi--homogeneous harmonic  polynomials in terms of the spaces  of so--called symplectic harmonics, which are defined by means of the scalar operators $\SE$ and $\SED$ and define irreducible $\mathrm{Sp}( p)$--representations.

\begin{definition}
If $a \geq b$ then 
\begin{equation*}
\SH_{a,b} := \mcH_{a,b}(\mR^{4p};\mC) \cap \mathrm{Ker} \, \SE\phantom{\dagger}
\end{equation*}
If $ a \leq b$ then 
\begin{equation*}
\SHD_{a,b} := \mcH_{a,b}(\mR^{4p};\mC) \cap \mathrm{Ker} \, \SED
\end{equation*}
\end{definition}

\begin{theorem}(see \cite{paper3})\\
For $a \geq b$, one has
$$
\mcH_{a,b}(\mR^{4p};\mC) = \SH_{a,b} \oplus \ \SED  \SH_{a+1,b-1} \oplus \ \SE^{\dagger 2}  \SH_{a+2,b-2} \oplus \cdots \oplus \ \SE^{\dagger b}  \SH_{a+b,0} 
$$
while for $a \leq b$, one has
$$
\mcH_{a,b}(\mR^{4p};\mC) = \SHD_{a,b} \oplus \ \SE \, \SHD_{a-1,b+1} \oplus \ \SE^2  \SHD_{a-2,b+2} \oplus \cdots \oplus \ \SE^{a}  \SHD_{0,b+a} 
$$
both decompositions being Fischer decompositions of the space of complex--valued harmonic bi--homogeneous  polynomials in terms of $\gsp_{2p}(\mC)$--irreducibles of complex--valued  (adjoint) symplectic harmonic  bi--homogeneous polynomials.
\label{theofischsymplharm}
\end{theorem}

Now we will show it to be unnecessary to distinguish between the cases $a \geq b$ and $a \leq b$ for obtaining a Fischer decomposition of the space $\mcH_{a,b}(\mR^{4p};\mC)$ in terms of symplectic spherical harmonics. Indeed, as we have seen, the Lie algebra $\gsl(2)$ can be realised in terms of the operators $\SE$ and $\SED$, and the spaces $\SH_{a,b}$ ($a \geq b$) and $\SHD_{a,b}$ ($a \leq b$) contain spherical harmonics of the indicated degree belonging to the kernel of these operators respectively. Algebraically speaking, this means that for each fixed pair $(a,b)$ of integers, these spaces can be seen as a highest weight vector (a lowest weight vector respectively) for a finite-dimensional irreducible $\gsl(2)$--module denoted by $\mV_{a-b}$ (see below). Note that, more precisely, there will be as many copies of this module as the dimension of the vector spaces $\SH_{a,b}$ ($\SHD_{a,b}$ respectively), since each polynomial inside these spaces generates an isomorphic copy of the aforementioned $\gsl(2)$-module, but we will always speak about `the' module generated by these vector spaces as a whole. 

\begin{lemma}
For fixed integers $a \geq b \geq 0$, the vector space $\SH_{a,b}$ generates, under the repeated action of the operator $\SED$, a model for the finite-dimensional irreducible $\gsl(2)$-module $\mV_{a - b}$ of dimension $(a - b + 1)$. The weight space decomposition then is given by 
$$ \mV_{a - b} = \bigoplus_{t=0}^{a - b}\SE^{\dagger t} \, \mcH_{a,b}^S$$
\end{lemma}

\noindent
{\bf Proof}\\
\noindent
As observed above the vector space $\mcH_{a,b}^S$ behaves as a highest weight vector, since it is annihilated by the operator $\mcE$ (which we regard as the positive root vector). The module then is generated by the repeated action of the complementary operator (the negative root vector $\mcE^\dagger$). In order to prove that this generates a module of the required dimension, we recall from \cite{paper3} that the mapping
$$
\mcE^{\dagger (a - b)} : \mcH_{a,b}^S \longrightarrow \mcH_{b,a}^{S\dagger} 
$$
defines an isomorphism, from which it follows that $\mcE^{\dagger(a - b+1)}\mcH_{a,b}^S = 0$. This identifies the space at the right--hand side as the lowest weight vector and hence fixes the dimension of $\mV_{a-b}$. \qed

\noindent
As a consequence we can now decompose the spaces $\mcH_{a,b}(\mR^m)$, with $a,b$ arbitrary, in terms of {\em either} only symplectic harmonics {\em or} only adjoint symplectic harmonics. Indeed, it suffices to note that each weight space in the module $\mV_{a-b}$ can be written using either the repeated action of the negative root vector $\SED$ on the highest weight vector $\SH_{a,b}$, or the repeated action of the positive root vector $\SE$ on the lowest weight vector $\SHD_{b,a}$. More precisely, the following identification holds.

\begin{lemma}
\label{identification}
For all integers $a \geq b \geq 0$ and $t=0,\ldots,a-b$, one has that
$$ 
\SE^{\dagger t} \, \SH_{a,b} \ \cong \ \SE^{a - b - t} \, \SHD_{b,a}
$$
\end{lemma}

\noindent
{\bf Proof}\\
\noindent
This immediately follows from standard $\gsl(2)$-relations and the fact that the symplectic (adjoint symplectic respectively) harmonics define the highest weight vector (the lowest weight vector respectively) for the $\gsl(2)$-module $\mV_{a-b}$. \qed

Taking into account Theorem \ref{theofischsymplharm} and Lemma \ref{identification} we obtain the following result.
\begin{proposition}
\label{decompharmsymplec}
For all bidegrees $(a,b)$ one has
$$
\mcH_{a,b}(\mR^{4p};\mC) = \bigoplus_{t=0}^a \, \SE^{\dagger (b-a+t)} \, \SH_{b+t,a-t}
$$
where from now on it is tacitly assumed that terms containing negative powers of $\SED$ are omitted.
\end{proposition}

Let us illustrate this result with an example. Consider the space $\mcH_{1,2}(\mR^{4p};\mC)$ which, according to Proposition \ref{decompharmsymplec}, can be decomposed as $\mcH_{1,2}(\mR^{4p};\mC) = \SED \, \SH_{2,1} \ \oplus  \ \SE^{\dagger 2} \, \SH_{3,0}$, or still, in view of Lemma \ref{identification}, as $\mcH_{1,2}(\mR^{4p};\mC) =  \SHD_{1,2} \ \oplus  \ \SE \, \SHD_{0,3}$. This exactly yields the decomposition of $\mcH_{1,2}(\mR^{4p};\mC)$ according to Theorem \ref{theofischsymplharm}.\\[-2mm]

In view of Proposition \ref{decompharmsymplec} we shall henceforth only use the spaces $\SH_{a,b}$ and the operator $\SED$, knowing that each time a similar, mirrored, result is valid involving the spaces $\SHD_{a,b}$ and the operator $\SE$. Along similar lines as in \cite{paper3}, the following proposition may then be proven.
\begin{proposition}
The space $\mcQ^{r,r}_{a,b}$ of q--monogenic functions with values in the symplectic cell $\mS^r_r$ may be decomposed into $\mathrm{Sp}( p)$--irreducibles as 
$$
\mcQ^{r,r}_{a,b} = \SE^{\dagger (b-a)}\mcS^{r}_{b,a} \, \oplus \, \SE^{\dagger (b-a+1)} \, \mcS^{r}_{b+1,a-1} \, \oplus \cdots \oplus \,  \SE^{\dagger b} \, \mcS^{r}_{b+a,0}
$$
\end{proposition}

If a q--monogenic function $F$ takes its values in the symplectic cell $\mS^r_s$ which is neither contained in Ker $P$ nor in Ker $Q$, then the function values can be shifted horizontally back and forth by means of the operators $P$ and $Q$, leading to the following result.

\begin{corollary}
The space $\mcQ^{r+2k,r}_{a,b}$ of q--monogenic functions with values in the symplectic cell $\mS^{r+2k}_r$ may be decomposed into $\mathrm{Sp}( p)$--irreducibles as
$$
\mcQ^{r+2k,r}_{a,b} = \SE^{\dagger (b-a)} \,  Q^k \, \mcS^{r}_{b,a} \, \oplus \, \SE^{\dagger (b-a+1)} \, Q^k \, \mcS^{r}_{b+1,a-1} \, \oplus \cdots \oplus \,  \SE^{\dagger b} \, Q^k \, \mcS^{r}_{b+a,0}
$$
\end{corollary}

Collecting all results found so far, the space $\mcP(\mR^{4p};\mS)$ of spinor--valued polynomials may be decomposed as
\begin{align*}
\mcP(\mR^{4p};\mS) &= \bigoplus_{k=0}^\infty \, \mcP_k(\mR^{4p};\mS) = \bigoplus_{k=0}^\infty \, \bigoplus_{\ell=0}^{\lfloor \frac{k}{2} \rfloor} \, r^{2\ell} \, \mcH_{k-2\ell}(\mR^{4p};\mS) \\
& = \bigoplus_{k=0}^\infty \, \bigoplus_{\ell=0}^{\infty} \, r^{2\ell} \, \mcH_{k}(\mR^{4p};\mS) 
 = \bigoplus_{k=0}^\infty \, \bigoplus_{\ell=0}^{\infty} \, r^{2\ell} \, \mcH_{k}(\mR^{4p};\mC) \otimes \mS\\
&=  \bigoplus_{\ell=0}^{\infty} \, \bigoplus_{a=0}^\infty \, \bigoplus_{b=0}^\infty \, |\uz|^{2\ell} \, \mcH_{a,b}(\mR^{4p};\mC)   \otimes \mS\\
&= \bigoplus_{\ell=0}^{\infty} \, \bigoplus_{a=0}^\infty \, \bigoplus_{b=0}^\infty \, |\uz|^{2\ell} \, \bigoplus_{t=0}^a \, \SE^{\dagger (b-a+t)} \, \SH_{b+t,a-t} \otimes  \bigoplus_{r=0}^{p}  \, \bigoplus_{j=0}^{p-r} \, Q^j \, \mS^{r}_r 
\end{align*}
or still
\begin{align*}
\mcP(\mR^{4p};\mS) &= \bigoplus_{\ell=0}^{\infty} \, \bigoplus_{a\geq b\geq 0}^\infty \, \bigoplus_{t=0}^{a-b} \, |\uz|^{2\ell} \,  \SE^{\dagger t} \, \SH_{a,b} \otimes  \bigoplus_{r=0}^{p}  \, \bigoplus_{j=0}^{p-r} \, Q^j \, \mS^{r}_r \\
&=  \bigoplus_{r=0}^{p}  \,   \bigoplus_{a\geq b\geq 0}^\infty \, \bigoplus_{t=0}^{a-b} \,  \bigoplus_{j=0}^{p-r} \, \bigoplus_{\ell=0}^{\infty} \, |\uz|^{2\ell} \, Q^j \,  \SE^{\dagger t} \, \SH_{a,b} \otimes    \mS^{r}_r 
\end{align*}

\noindent So we now need to obtain a Fischer decomposition of the spaces $\SH_{a,b,r} := \SH_{a,b} \otimes \mS^r_r$, consisting of bi--homogeneous harmonic polynomials in the kernel of the operators $\SE$ and $P$, in terms of $\mathrm{Sp}( p)$--irreducibles.\\[-2mm]

Let us give a small flavour of this Fischer decomposition by means of the following simple example. Take $p=2$ and consider the first order polynomial $\mathrm{P}_{1,0} = z_2 \, \gfd_1 \, I  \in \SH_{1,0,1}$, which we expect (see below) to decompose as
\begin{equation}
\label{easyexample}
\mathrm{P}_{1,0} = z_2 \, \gfd_1 \, I = S_{1,0}^1 + \uz \, S_{0,0}^2 + (\uzJ + A \uz \, Q) \, S_{0,0}^0
\end{equation}
Expressing that $\mathrm{P}_{1,0}$ belongs to Ker $P$ determines $A = -\onehalf$; expressing its q--monogenicity determines $S_{0,0}^0= \frac{1}{6} \, I$ and $S_{0,0}^2 = \frac{1}{4} (- \gfd_1 \gfd_2 + \gfd_3 \gfd_4) I$, and hence
\begin{equation}
\label{decompeasyexample}
\mathrm{P}_{1,0} = z_2 \, \gfd_1 \, I = \onehalf (z_1 \gfd_2 + z_2 \gfd_1) + \uz \, \frac{1}{4} (- \gfd_1 \gfd_2 + \gfd_3 \gfd_4) I + (\uzJ -\onehalf \uz \, Q)  \frac{1}{6} \, I
\end{equation}

It remains to explain how the predicted form (\ref{easyexample}) of the decomposition was obtained.
If the given polynomial is in $\SH_{a,b,r}$ (with $a>b$) then the first component in its Fischer decomposition will turn out to be the so--called Cartan piece $\mcS_{a,b}^r$, see the forthcoming paper \cite{paperdirect}; let us call it the {\em zero--letter--word component}, since it is a polynomial space multiplied by a constant embedding factor, and denote it by $\mcS_{a,b}^{r,0}$.\\[-2mm]

The next term in the decomposition is called the {\em one--letter--word component} and is the direct sum of the projections of the spaces $\uz \, \mcS_{a-1,b}^{r+1}$, $\uzJd \, \mcS_{a,b-1}^{r+1}$, $\uzd \, \mcS_{a,b-1}^{r-1}$, and $\uzJ \, \mcS_{a-1,b}^{r-1}$ on Ker $\Delta$ (this projection being trivial), on Ker $\SE$ and on Ker $P$. Let us denote by  $T_{a,b}^r$ a polynomial in one of these spaces. It is clear that $T_{a,b}^r$ takes values either in $\mS^r_r$ or in $\mS^r_r \oplus \mS_{r-2}^r$, whence $PT_{a,b}^r$ takes values in $\mS_{r-2}^{r-2}$ and $P^2 T_{a,b}^r = 0$. The projection of $T_{a,b}^r$ on Ker $P$ thus is given by
$$
\pi_{\mathrm{Ker} P}[ T_{a,b}^r ] = \left( {\bf 1} - \frac{1}{p-r+2} QP   \right) T_{a,b}^r
$$
which explicitly yields
\begin{align*}
\pi_{\mathrm{Ker} P}[\uz \, S_{a-1,b}^{r+1}] &= \uz \, S_{a-1,b}^{r+1}\\
\pi_{\mathrm{Ker} P}[\uzJd \, S_{a,b-1}^{r+1}] &= \uzJd \, S_{a,b-1}^{r+1}\\
\pi_{\mathrm{Ker} P}[\uzd \, S_{a,b-1}^{r-1}] &= \left( \uzd + \frac{1}{p-r+2} Q \uzJd   \right) S_{a,b-1}^{r-1}\\
\pi_{\mathrm{Ker} P}[\uzJ \, S_{a-1,b}^{r-1}] &= \left( \uzJ - \frac{1}{p-r+2} Q \uz  \right) S_{a-1,b}^{r-1}
\end{align*}
Next, the projection on Ker $\SE$ of a polynomial $T^r_{a,b}$ for which $\SE^2 T^r_{a,b} = 0$ is given by
$$
\pi_{\mathrm{Ker} \SE}[ T_{a,b}^r ] = \left( {\bf 1} - \frac{1}{a-b+2} \SED \SE   \right) T_{a,b}^r
$$
leading to
\begin{align*}
\pi_{\mathrm{Ker} \SE}[\uz \, S_{a-1,b}^{r+1}] &= \uz \, S_{a-1,b}^{r+1}\\
\pi_{\mathrm{Ker} \SE}[\uzJd \, S_{a,b-1}^{r+1}] &= \left( \uzJd - \frac{1}{a-b+2} \SED \uz \right) S_{a,b-1}^{r+1}\\
\pi_{\mathrm{Ker} \SE}[\uzd \, S_{a,b-1}^{r-1}] &= \left( \uzd + \frac{1}{a-b+2} \SED \uzJ   \right) S_{a,b-1}^{r-1}\\
\pi_{\mathrm{Ker} \SE}[\uzJ \, S_{a-1,b}^{r-1}] &= \uzJ \, S_{a-1,b}^{r-1}
\end{align*}
The one--letter--word component in the Fischer decomposition thus will be a direct sum of the following polynomial spaces:
\begin{align*}
\mcS_{a,b}^{r,1} &:= \pi_{\mathrm{Ker} \SE} \, \pi_{\mathrm{Ker} P} \, [\uz \, \mcS_{a-1,b}^{r+1}] = \uz \, \mcS_{a-1,b}^{r+1}\\
\mcS_{a,b}^{r,2} &:= \pi_{\mathrm{Ker} \SE} \, \pi_{\mathrm{Ker} P} \, [\uzJd \, \mcS_{a,b-1}^{r+1}] = \left( \uzJd - \frac{1}{a-b+2} \SED \uz \right) \mcS_{a,b-1}^{r+1}\\
\mcS_{a,b}^{r,3} &:= \pi_{\mathrm{Ker} \SE} \, \pi_{\mathrm{Ker} P} \, [\uzd \, \mcS_{a,b-1}^{r-1}]\\  
&\phantom{:} = \left( \uzd + \frac{1}{a-b+2} \SED \uzJ  +
\frac{1}{p-r+2} Q \, \uzJd - \frac{1}{(p-r+2) (a-b+2)} Q \, \SED \, \uz \right) \mcS_{a,b-1}^{r-1}\\
\mcS_{a,b}^{r,4} &:= \pi_{\mathrm{Ker} \SE} \, \pi_{\mathrm{Ker} P} \, [\uzJ \, \mcS_{a-1,b}^{r-1}] = \left( \uzJ - \frac{1}{p-r+2} Q \, \uz  \right) \mcS_{a-1,b}^{r-1}
\end{align*}

\noindent
Now we turn our attention to the {\em two--letter--word component} in the  Fischer decomposition, which, a priori, is the direct sum of the following combinations: $\uz \uzd \mcS^r_{a-1,b-1}$, $\uzJ \uzJd \mcS^r_{a-1,b-1}$, $\uz \uzJd \mcS^{r+2}_{a-1,b-1}$, $\uzJ \uzd \mcS^{r-2}_{a-1,b-1}$, $\uz \uzJ \mcS^r_{a-2,b}$ and $\uzd \uzJd \mcS^r_{a,b-2}$, or more precisely, of their images under the composition $\pi = \pi_{\mathrm{Ker} \SE}  \circ\pi_{\mathrm{Ker} P} \circ  \pi_{\mathrm{Ker} \Delta}$ of the projections on the kernels of $\Delta$, $P$ and $\SE$. For these projections the following formulae are valid:
\begin{align*}
\pi_{\mathrm{Ker} \Delta}[T^r_{a,b}] &= T^r_{a,b} - \frac{1}{4(2p+a+b-2)} r^2  \Delta T^r_{a,b} + \frac{1}{32(2p+a+b-2)(2p+a+b-3)} r^4 \Delta^2 T^r_{a,b}\\
\pi_{\mathrm{Ker} P}[T^r_{a,b}] &= T^r_{a,b} - \frac{1}{p-r+2} Q P T^r_{a,b} + \frac{1}{2(p-r+3)(p-r+2)} Q^2  P^2 T^r_{a,b}\\
\pi_{\mathrm{Ker} \SE}[T^r_{a,b}] &= T^r_{a,b} - \frac{1}{a-b+2} \SED \SE T^r_{a,b} + \frac{1}{2(a-b+3)(a-b+2)} \SE^{\dagger 2} \SE^2 T^r_{a,b}
\end{align*}

\noindent Straightforward calculations then show the two--letter--word component to be the direct sum of the following polynomial spaces:
\begin{align*}
\mcS_{a,b}^{r,5} &:= \pi[  \uz \uzd \mcS^r_{a-1,b-1} ]\\
 &\phantom{:} = \left(   \uz \uzd   + \frac{1}{a-b+2} \SED \uz \uzJ + \frac{1}{p-r+2} Q \uz \uzJd - \frac{2p+b-r-1}{2p+a+b-2} |\uz|^2  \right)  \mcS^r_{a-1,b-1} \\
\mcS_{a,b}^{r,6} &:= \pi[  \uzJ \uzJd \mcS^r_{a-1,b-1}  ]\\ 
& \phantom{:} =  \left(   \uzJ \uzJd  - \frac{1}{a-b+2} \SED \uzJ \uz - \frac{1}{p-r+2} Q \uz \uzJd  - \frac{b+r-1}{2p+a+b-2} |\uz|^2 \right) \mcS^r_{a-1,b-1} \\
\mcS_{a,b}^{r,7} &:= \pi[  \uz \uzJd \mcS^{r+2}_{a-1,b-1} ] =  \left(    \uz \uzJd     \right) \mcS^{r+2}_{a-1,b-1} 
\end{align*}
\begin{align*}
\mcS_{a,b}^{r,8} &:= \pi[  \uzJ \uzd \mcS^{r-2}_{a-1,b-1} ]\\ 
 & \phantom{:}=  \left (   \uzJ \uzd  - \frac{1}{p-r+2} Q \uz \uzd + \frac{1}{p-r+2} Q \uzJ \uzJd 
- \frac{1}{(p-r+3)(p-r+2)} Q^2 \uz \uzJd \right .\\
 & \hspace*{15mm} + \left .\frac{1}{(p-r+3)(p-r+2)(a-b+2)}  \SED Q^2  \uzJ \uz \right ) \, \mcS^{r-2}_{a-1,b-1} \\[2mm]
\mcS_{a,b}^{r,9} &:= \pi[  \uz \uzJ \mcS^r_{a-2,b} ] =  \left(  \uz \uzJ   \right) \mcS^r_{a-2,b} \\
\mcS_{a,b}^{r,10} &:= \pi[  \uzd \uzJd \mcS^r_{a,b-2} ]\\ 
& \phantom{:}=  \left(  \uzd \uzJd  - \frac{1}{a-b+2} \SED (\uzd \uz - \uzJ \uzJd) - \frac{1}{(a-b+3)(a-b+2)}  \SE^{\dagger 2} \uzJ \uz \right) \mcS^r_{a,b-2}
\end{align*}

As to the three--letter--word component, we find the following terms:
\begin{align*}
\mcS_{a,b}^{r,11} &:= \pi \left[ \uz \uzd \uzJ \mcS_{a-2,b-1}^{r-1} \right]\\
 & \phantom{:}= \left( \uz \uzd \uzJ - \frac{2p+b+1-r}{2p+a+b-2} \, |\uz|^2 \uzJ
- \frac{1}{p+2-r} \, Q \uz \uzJ \uzJd \right . \\
& \hspace*{17mm} + \left . \frac{2p+b+1-r}{(p+2-r)(2p+a+b-2)} \, Q  | \uz |^2 \uz \right) \mcS_{a-2,b-1}^{r-1} \\
\mcS_{a,b}^{r,12} &:= \pi \left[ \uz \uzd \uzJd \mcS_{a-1,b-2}^{r+1} \right]\\ 
& \phantom{:}= \left(  \uz \uzd \uzJd - \frac{2p+b-2-r}{2p+a+b-2} \,  | \uz |^2 \uzJd 
- \frac{1}{a-b+2} \, \SED \uz \uzJd \uzJ  \right . \\
& \hspace*{18mm} + \left . \frac{2p+b-2-r}{(a-b+2)(2p+a+b-2)} \, \SED  | \uz |^2 \uz   \right) \mcS_{a-1,b-2}^{r+1}  \\
\mcS_{a,b}^{r,13} &:= \pi \left[ \uz \uzJ \uzJd \mcS_{a-2,b-1}^{r+1} \right]= \left( \uz \uzJ \uzJd - \frac{b+r-1}{2p+a+b-2} \, |\uz|^2 \uz  \right) \mcS_{a-2,b-1}^{r+1}
\end{align*}
\begin{align*}
\mcS_{a,b}^{r,14} &:= \pi \left[ \uzd \uzJ \uzJd \mcS_{a-1,b-2}^{r-1} \right]\\ 
& \phantom{:}= \left ( \uzd \uzJ \uzJd - \frac{b+r-4}{2p+a+b-2} \, |\uz|^2 \uzd 
- \frac{1}{a-b+2} \, \SED \uzd \uzJ \uz + \frac{1}{p-r+2} \, Q \uz \uzd \uzJd \right .\\ 
& \hspace*{20mm} - \frac{b+r-4}{(2p+a+b-2)(a-b+2)} \, \SED |\uz|^2 \uzJ - \frac{b+r-4}{(2p+a+b-2)(p-r+2)} \, Q |\uz|^2 \uzJd \\ 
& \hspace*{20mm}  - \frac{1}{(p-r+2)(a-b+2)} \SED Q \uz \uzJd \uzJ \\
& \hspace*{20mm} + \left . \frac{b+r-4}{(2p+a+b-2)(p-r+2)(a-b+2)} \, \SED Q |\uz|^2 \uz \right ) \, \mcS_{a-1,b-2}^{r-1}
\end{align*}

Finally the four--letter--word component turns out to be
\begin{align*}
\mcS_{a,b}^{r,15} &:= \pi \left[ \uz \uzd \uzJ \uzJd \mcS_{a-2,b-2}^r \right] = \left(\uz \uzd \uzJ \uzJd - \frac{1}{4(2p+a+b-2)} \, |\uz|^2 \, \mcT_{1,1} + \omega |\uz|^4\right)  \mcS_{a-2,b-2}^r
\end{align*}
with
$$
\mcT_{1,1} = 4(b+r-4) \uz \uzd + 4(2p+b-r)\uzJ\uzJd
- 4 \SED \uz \uzJ  - 8 Q \uz \uzJd
$$
and
$$
\omega = \frac{2pb + b^2 - 5b - a + 2pr - 2r - r^2 - 8p+6}{(2p+a+b-2)(2p+a+b-3)} 
$$

Now we claim to have found all components in the Fischer decomposition, with respect to the Sp$( p)$--action, of the space $\SH_{a,b,r} = \SH_{a,b} \, \otimes \, \mS^r_r$, in other words we do not have to consider five--or--more--letter word components. Indeed, it is easily observed that a five--letter--word immediately reduces to a three--letter--word multiplied by $|\uz|^2$, which precisely is one of the operators appearing in the Howe dual pair which will be used to assemble isomorphic Sp$( p)$--irreducibles. This observation moreover is supported by an abstract reasoning  based on group representation theoretical arguments which will be developed in detail in \cite{paperdirect}. In this way the following result is obtained.

\begin{proposition}
With respect to the $\mathrm{Sp}( p)$--action, the space $\SH_{a,b,r}$ of bi--homogeneous symplectic harmonic  polynomials with values in the symplectic cell $\mS^r_r$, can  be decomposed into irreducibles as
$$
\SH_{a,b,r}  =  \bigoplus_{\alpha=0}^{15} \, \mcS_{a,b}^{r,\alpha}
$$
\end{proposition}

Returning to the decomposition of $\mcP(\mR^{4p};\mS)$ in terms of spaces of symplectic harmonics, we finally obtain its decomposition in terms of spaces of $\gosp(4|2)$--monogenics.

\begin{theorem}
\label{finaltheo}
With respect to the $\mathrm{Sp}( p)$--action, the space $\mcP(\mR^{4p};\mS)$ of spinor--valued polynomials, may be decomposed into irreducibles as 
$$
\mcP(\mR^{4p};\mS) = \bigoplus_{r=0}^{p}  \,   \bigoplus_{a\geq b\geq 0}^\infty \, \bigoplus_{\alpha=0}^{15} \, \bigoplus_{t=0}^{a-b} \,  \bigoplus_{j=0}^{p-r} \, \bigoplus_{\ell=0}^{\infty} \,  |\uz|^{2\ell} \, Q^j \,  \SE^{\dagger t} \,  \mcS_{a,b}^{r,\alpha}
$$
\end{theorem}


\section*{Acknowledgements}
D.\ Eelbode was supported by the Fund for Scientic Research-Flanders (FWO-V), on the project "Construction of algebra realisations using Dirac-operators", grant G.0116.13N. \\
R.\ L\'avi\v{c}ka and V.\ Sou\v{c}ek gratefully acknowledge support by the Czech grant GA CR G201/12/028. \\
This paper was finalized during a scientific stay of F.\ Brackx, H.\ De Schepper and D.\ Eelbode at the Mathematical Institute of the Charles Univerity in Prague; they heartily want to thank their co-authors for their hospitality at that, and many other occasions.


\end{document}